\documentclass[centertags,10pt,twoside,onecolumn]{amsart}
\usepackage[a4paper, textwidth=130mm, textheight=230mm]{geometry}
\calclayout

\usepackage{tikz}
\usetikzlibrary{arrows.meta}
\usepackage{subcaption}

\usepackage[flushmargin]{footmisc}

\usepackage{amsmath, amsthm, amssymb, enumitem, microtype, mlmodern}

\usepackage{xcolor}
\usepackage[colorlinks]{hyperref}
\definecolor{linkcolor}{HTML}{800006}
\definecolor{citecolor}{HTML}{2E7E2A}
\definecolor{filecolor}{HTML}{131877}
\definecolor{urlcolor} {HTML}{690067}
\definecolor{menucolor}{HTML}{727500}
\definecolor{runcolor} {HTML}{137776}
\hypersetup{
    linkcolor=linkcolor,
    citecolor=citecolor,
    filecolor=filecolor,
    urlcolor=urlcolor,
    menucolor=menucolor,
    runcolor=runcolor
}

\theoremstyle{plain}
\newtheorem{theo}{Theorem}[section]
\newtheorem*{theo*}{Theorem}
\newtheorem{conj}[theo]{Conjecture}
\newtheorem{lemm}[theo]{Lemma}
\newtheorem{prop}[theo]{Proposition}
\newtheorem*{prop*}{Proposition}
\newtheorem{coro}[theo]{Corollary}
\newtheorem{claim}[theo]{Claim}

\theoremstyle{definition}
\newtheorem{defi}[theo]{Definition}
\newtheorem{exam}[theo]{Example}

\usepackage[style=alphabetic, isbn=false, doi=false]{biblatex}
\newbibmacro{string+doi}[1]{%
	\iffieldundef{doi}{\iffieldundef{url}{#1}{\href{\thefield{url}}{#1}}}{\href{http://dx.doi.org/\thefield{doi}}{#1}}}	
\DeclareFieldFormat[article]{title}{\usebibmacro{string+doi}{\mkbibquote{#1}}}
\DeclareFieldFormat[book]{title}{\usebibmacro{string+doi}{\mkbibemph{#1}}}

\newcommand{\vocab}[1]{\textbf{#1}}
\newcommand{\integers}{\mathbb{Z}}
\newcommand{\positiveintegers}{\mathbb{Z}_{>0}}
\newcommand{\mtp}[1]{^{(#1)}} \newcommand{\graph}{\varGamma}
\newcommand{\vertices}{V}
\newcommand{\edges}{E}
\newcommand{\powerset}[1]{\mathcal{P}(#1)}
\newcommand{\coefficientring}{R}
\newcommand{\laurentpolynomialring}{\mathcal{L}}
\newcommand{\laurentmonomial}[2]{\ell(#1, #2)}
\DeclareMathOperator{\normalizedweight}{nwt}
\DeclareMathOperator{\weight}{wt}
\newcommand{\sumacy}[1]{\sum_{#1}^{\mathsf{acy}}}
\newcommand{\graphLPalgebra}[1]{\mathcal{A}_{#1}}
\newcommand{\monomial}[2]{m(#1, #2)}
\newcommand{\coefficient}[2]{c(#1, #2)}
\newcommand{\ppath}[3]{\ensuremath{#1 \to_{#2} #3}}

\title
[Cluster monomials in graph LP algebras]
{Cluster monomials in graph Laurent phenomenon algebras}

\author[G. Z. Dantas e Moura]{Guilherme Zeus Dantas e Moura}
\address{Guilherme Zeus Dantas e Moura, \normalfont University of Waterloo,
Department of Combinatorics and Optimization,
200 University Avenue West,
Waterloo, ON N2L 3G1 (Canada)}
\email{zeus.dantasemoura@uwaterloo.ca}
\email{zeus@guilhermezeus.com}
\urladdr{https://www.guilhermezeus.com}

\author[R. C. Telekicherla Kandalam]{Ramanuja Charyulu Telekicherla Kandalam}
\address{Ramanuja Charyulu Telekicherla Kandalam, \normalfont University of Minnesota,
School of Mathematics,
206 Church Street SE,
Minneapolis, MN 55455 (USA)}
\email{telek002@umn.edu}

\author[D. Woodruff]{Dora Woodruff}
\address{Dora Woodruff, \normalfont Massachusetts Institute of Technology,
Department of Mathematics,
77 Massachusetts Avenue,
Cambridge, MA 02139 (USA)}
\email{dorawood@mit.edu}

\thanks{\textit{Acknowledgements.}
This project was partially supported by RTG grant NSF/DMS-1745638.
It was supervised as part of the University of Minnesota School of Mathematics Summer 2023 REU program.
The first author was partially supported by Haverford College's KINSC Summer Scholar funding and by NSERC grant RGPIN-2021-02568.
Part of this work was completed while the first author was affiliated with Haverford College and the third author was affiliated with Harvard University.}

  

\begin{document}

\begin{abstract}
    Laurent phenomenon algebras,
    first introduced by Lam and Pylyavskyy,
    are a generalization of cluster algebras that still possess many salient features of cluster algebras.
    Graph Laurent phenomenon algebras,
    defined by Lam and Pylyavskyy,
    are a subclass of Laurent phenomenon algebras
    whose structure is given by the data of a directed graph.
    In this paper,
    we prove that the cluster monomials of a graph Laurent phenomenon algebra form a linear basis,
    as conjectured by Lam and Pylyavskyy
    and analogous to a result for cluster algebras by Caldero and Keller.
    We also prove that, if the graph is a bidirected tree,
    the coefficients of the expansion of any monomial in terms of cluster monomials are nonnegative.
\end{abstract}

\maketitle

\section{Introduction} \label{section:introduction}

Fomin and Zelevinsky~\cite{clusteralgebrasi} introduced cluster algebras, which have since been recognized as rather ubiquitous in mathematics.
Originally introduced as a combinatorial model for total positivity, cluster algebras are relevant to Teichmüller theory~\cite{cluster-teich}, triangulated surfaces~\cite{cluster-triangulated-surface}, Lie theory, representation theory of quivers, and Poisson geometry~\cite{cluster-poisson}, to name some examples.

Broadly, in a cluster algebra, a set of generators called the \textit{cluster variables} are organized into sets called \textit{clusters}. A \textit{seed} consists of a cluster and an \textit{exchange polynomial} for each cluster variable in the cluster. These exchange polynomials are always \emph{binomials} in the variables of the cluster and define a procedure for mutating a seed into a different seed. 

The Laurent phenomenon is a remarkable property of cluster algebras. It states that any cluster variable can be expressed as a \textit{Laurent polynomial} when written as a rational function in any cluster. Motivated by this property, Lam and Pylyavskyy~\cite{LP} defined \emph{Laurent phenomenon algebras}. Laurent phenomenon algebras are generalizations of cluster algebras in which the exchange polynomials are no longer required to be binomials. Lam and Pylyavskyy~\cite{LP} show that in their setting, the Laurent phenomenon still holds. 

Lam and Pylyavskyy~\cite{linearLP} also studied a specific class of Laurent phenomenon algebras called \textit{graph Laurent phenomenon algebras}.
In a graph Laurent phenomenon algebra, the exchange polynomials are linear polynomials determined by some directed graph \(\graph\).
Lam and Pylyavskyy~\cite{linearLP} provided an explicit combinatorial description of the cluster variables and clusters of a graph Laurent Phenomenon algebra \(\graphLPalgebra{\graph}\) (see Subsection~\ref{subsection:graphLPalgebraandclusters}).

The main result of this paper is Theorem~\ref{theorem:fromlinearLP}, which confirms a conjecture by \cite[Conjecture~7.3a]{linearLP} that the cluster monomials form a linear basis of \(\graphLPalgebra{\graph}\).

\begin{theo}\label{theorem:fromlinearLP}
    Let \(\graph\) be a graph.
    The cluster monomials of the graph LP algebra \(\graphLPalgebra{\graph}\) with coefficient ring \(\coefficientring\) form a linear basis of \(\graphLPalgebra{\graph}\) over \(\coefficientring\).
\end{theo}

Theorem~\ref{theorem:fromlinearLP} is motivated by an analogous fact about cluster algebras: Caldero and Keller~\cite{basis} proved that cluster monomials form a linear basis in finite type cluster algebras.

The proof of Theorem~\ref{theorem:fromlinearLP} is split into two parts.
First, we establish the linear independence of cluster monomials in a graph LP algebra, as stated in Theorem~\ref{theorem:linearindependence}.
Second, we demonstrate that cluster monomials in a graph LP algebra form an \(\coefficientring\)-linear spanning set, as stated in Theorem~\ref{theorem:spanningset}. 

Lam and Pylyavskyy~\cite{linearLP} also propose Conjecture~\ref{conjecture:positivity},
a stronger version of Theorem~\ref{theorem:fromlinearLP} that states not only that any monomial of \(\graphLPalgebra{\graph}\) is a linear combination of cluster monomials over \(\coefficientring\),
but also that the coefficients are nonnegative.
This stronger conjecture is also motivated by analogous results in the case of cluster algebras, by Lee and Schiffler~\cite{cluster-positivity-LS}.

\begin{conj}[{\cite[Conjecture 7.3b]{linearLP}}] \label{conjecture:positivity}
    Each monomial of \(\graphLPalgebra{\graph}\) is a linear combination of cluster monomials over \(\coefficientring\) with nonnegative coefficients.
\end{conj}

We prove Conjecture~\ref{conjecture:positivity} for the special case when the graph \(\graph\) is a bidirected tree, as stated in Theorem~\ref{theorem:tree}.

\begin{theo} \label{theorem:tree}
    Assume that \(\graph\) is a bidirected tree. Then, each monomial of \(\graphLPalgebra{\graph}\) is a linear combination of cluster monomials over \(\coefficientring\) with nonnegative coefficients.
\end{theo}

The article is organized as follows.
In Section \ref{section:preliminaries}, we establish notation and recall definitions and results from \cite{linearLP}.
In Section \ref{section:linearindependence}, we prove that cluster monomials in a graph LP algebra are linearly independent over \(\coefficientring\).
In Section \ref{section:spanningset}, we prove that cluster monomials in a graph LP algebra form an \(\coefficientring\)-linear spanning set. Finally, in Section \ref{section:trees}, we prove the nonnegativity result for trees. 

\section{Preliminaries and notational conventions} \label{section:preliminaries}

\subsection{Multisets} \label{subsection:multisets}

A \vocab{multiset} \(S\) is a sequence of sets \(S\mtp{1} \supset S\mtp{2} \supset \cdots\) indexed by the positive integers.
We say that \(x\) has \vocab{multiplicity} \(m\) in \(S\) if \(m+1\) is the smallest positive integer such that \(x \notin S\mtp{m+1}\).
In this case, if \(m > 0\), it follows that \(x \in S\mtp{m}\).
Therefore, \(S\mtp{i}\) is the set of elements with multiplicity at least \(i\) in \(S\).
A set \(S\) is naturally identified with the multiset given by the sequence \(S \supset \varnothing \supset \varnothing \supset \cdots\), that is, we interpret a set as a multiset with all multiplicities equal to \(1\).

Given multisets \(S\) and \(T\), we say that \(S\) and \(T\) are \vocab{disjoint} if \(S\mtp{1}\) and \(T\mtp{1}\) are disjoint.
We say that \(S\) is \vocab{contained} in \(T\), denoted by \(S \subset T\), if \(S\mtp{i} \subset T\mtp{i}\) for all \(i\).
The \vocab{sum} of two multisets \(S\) and \(T\) is the multiset \(S + T\) with the multiplicity of each element \(x\) being the sum of its multiplicities in \(S\) and \(T\).
If \(S\) is contained in \(T\), the \vocab{subtraction} of \(S\) from \(T\) is the multiset \(T - S\) with the multiplicity of each element \(x\) being the subtraction of its multiplicity in \(S\) from its multiplicity in \(T\).

We may define a multiset by simply listing its elements, with the multiplicity of each element being the number of times it appears in the list.
For example, \(S = \{1, 1, 2, 3\}\) is the multiset with elements \(1, 2, 3\) and multiplicities \(2, 1, 1\), respectively, in which case we have \(S\mtp{1} = \{1, 2, 3\}\), \(S\mtp{2} = \{1\}\), and \(S\mtp{i} = \varnothing\) for all \(i \geq 3\).

\subsection{Directed graph} \label{subsection:directedgraph}

Let \(\graph\) be a \vocab{directed graph} with vertex set \(\vertices\) and edge set \(\edges\).
We maintain this notation throughout the document.
Examples of directed graphs are given in Figure~\ref{fig:directedgraphs}.

\begin{figure}[htbp]
    \centering
    \begin{subfigure}[b]{0.47\textwidth}
        \centering
        \begin{tikzpicture}
            \node[shape=circle, draw=black] (2) at (0,0){2};
            \node[shape=circle, draw=black] (3) at (-1,-1){3};
            \node[shape=circle, draw=black] (1) at (-1,1){1};
            \node[shape=circle, draw=black] (4) at (-2,0){4};
            \path[draw,thick]
                (1) edge node {} (2)
                (3) edge node {} (1)
                (3) edge node {} (4)
                (3) edge node {} (2)
                (4) edge node {} (3)
                (4) edge node {} (1);
        \end{tikzpicture}
        \caption{A directed graph with four vertices, \(1\), \(2\), \(3\), and \(4\), and ten directed edges, \(12\), \(13\), \(14\), \(21\), \(23\), \(31\), \(32\), \(34\), \(41\), and \(43\).}
        \label{fig:directedgraphs:1}
    \end{subfigure}
    \hfill
    \begin{subfigure}[b]{0.47\textwidth}
        \centering
        \begin{tikzpicture}
            \node[shape=circle, draw=black] (2) at (0,0){2};
            \node[shape=circle, draw=black] (3) at (-1,-1){3};
            \node[shape=circle, draw=black] (1) at (-1,1){1};
            \node[shape=circle, draw=black] (4) at (-2,0){4};
            \path[draw,thick, ->] (1) edge node {} (2);
            \path[draw,thick, ->] (2) edge node {} (3);
            \path[draw,thick, ->] (3) edge node {} (4);
            \path[draw,thick, ->] (4) edge node {} (1);
            \path[draw,thick] (1) edge node {} (3);    
        \end{tikzpicture}
        \caption{A directed graph with four vertices, \(1\), \(2\), \(3\), and \(4\), and six directed edges, \(12\), \(13\), \(23\), \(31\), \(34\), and \(41\).}
        \label{fig:directedgraphs:2}
    \end{subfigure}
    \caption{Two directed graphs.  An undirected edge represents two directed edges, one in each direction. An undirected graph will mean a graph in which all edges are undirected.}
    \label{fig:directedgraphs}
\end{figure}

In general, whenever we refer to an undirected graph, we mean that every edge is bidirected. 

\subsection{Nested collections} \label{subsection:nestedcollections}

A subset \(I \subset \vertices\) is \vocab{strongly connected} if the induced subgraph on \(I\) is strongly connected, that is, if for all vertices \(v, u \in I\), there is some directed path contained in \(I\) from \(v\) to \(u\).

For example,
taking \(\graph\) to be the directed graph in Figure~\ref{fig:directedgraphs:1},
the strongly connected sets of vertices are \(\varnothing\), \(\{1\}\), \(\{2\}\), \(\{3\}\), \(\{4\}\), \(\{1, 2\}\), \(\{1, 3\}\), \(\{1, 4\}\), \(\{2, 3\}\), \(\{3, 4\}\), \(\{1, 2, 3\}\), \(\{1, 2, 4\}\), \(\{1, 3, 4\}\), \(\{2, 3, 4\}\), and \(\{1, 2, 3, 4\}\).
Taking \(\graph\) to be the directed graph in Figure~\ref{fig:directedgraphs:2},
the strongly connected sets of vertices are \(\varnothing\), \(\{1\}\), \(\{2\}\), \(\{3\}\), \(\{4\}\), \(\{1, 3\}\), \(\{1, 2, 3\}\), \(\{1, 3, 4\}\), and \(\{1, 2, 3, 4\}\).

Note that any subset \(U \subset \vertices\) is uniquely partitioned into maximal strongly connected subsets, called \vocab{strongly connected components} of \(U\).

A (multi)set \(\mathcal{N}\) with elements from \(\powerset{\vertices}\), where \(\powerset{\vertices}\) denotes the powerset of \(\vertices\),
is a \vocab{nested (multi)collection} if
\begin{enumerate}[label = \textup{(N\arabic*)}, ref = \textup{(N\arabic*)}]
    \item \label{item:nestedcollection:disjoint}
    for every pair \(I, J \in \mathcal{N}\),
    either \(I \subset J\), \(J \subset I\), or \(I \cap J = \varnothing\), and
    \item \label{item:nestedcollection:stronglyconnected}
    for any \(\mathcal{R} \subset \mathcal{N}\) such that \(I \cap J = \varnothing\) for all distinct \(I, J \in \mathcal{R}\),
    each \(I \in \mathcal{R}\) is a strongly connected component of the subgraph induced by \(\bigcup_{J \in \mathcal{R}} J\).
\end{enumerate}

Note that \ref{item:nestedcollection:stronglyconnected}, applied to \(\mathcal{R} = \{I\}\), implies that each \(I \in \mathcal{N}\) is strongly connected.
Note that a multiset \(\mathcal{N}\) is a nested multicollection if and only if the set \(\mathcal{N}\mtp{1}\) is a nested collection.
We usually disregard the empty set when considering nested collections.

For example,
consider \(\graph\) to be the directed graph in Figure~\ref{fig:directedgraphs:1}.
The set
\begin{equation*}
    \{
        \{2\},
        \{4\},
        \{2, 3, 4\},
        \{1, 2, 3, 4\}
    \}
\end{equation*} 
is a nested collection.
The set
\begin{equation*}
    \{
        \{1\},
        \{1, 3\},
        \{1, 4\}
    \}
\end{equation*}
is not a nested collection because \(\{1, 3\}\) and \(\{1, 4\}\) do not satisfy \ref{item:nestedcollection:disjoint}.
The set
\begin{equation*}
    \{
        \{2, 4\},
        \{2, 3, 4\},
        \{1, 2, 3, 4\}
    \}
\end{equation*}
is not a nested collection because \(\mathcal{R} = \{ \{2, 4\} \}\) does not satisfy \ref{item:nestedcollection:stronglyconnected}, since the strongly connected components of \(\{2, 4\}\) are \(\{2\}\) and \(\{4\}\).
The set
\begin{equation*}
    \{
        \{1\},
        \{3\},
        \{1, 2, 3\},
        \{1, 2, 3, 4\}
    \}
\end{equation*}
is not a nested collection because \(\mathcal{R} = \{ \{1\}, \{3\} \}\) does not satisfy \ref{item:nestedcollection:stronglyconnected}, since the only strongly connected component of \(\{1, 3\}\) is \(\{1, 3\}\).

\begin{lemm} \label{lemma:multinested:unique}
    The map \(\mathcal{N} \mapsto \sum_{I \in \mathcal{N}} I\) is a bijection between the set of nested multicollections \(\mathcal{N}\) and the set of multisets \(T\) with elements from \(\vertices\).
\end{lemm}

The proof of Lemma~\ref{lemma:multinested:unique} is split into two parts: surjectivity and injectivity.

\begin{proof}[Proof (surjectivity)]
    Let \(T\) be a multiset with elements from \(\vertices\).
    We construct a nested multicollection \(\mathcal{N}\) such that \(T = \sum_{I \in \mathcal{N}} I\).

    Recall that \(T\mtp{i}\) denotes the set of elements of \(T\) with multiplicity at least \(i\).
    Note that \(T\mtp{1} \supset T\mtp{2} \supset \cdots\).
    Let \(\mathcal{N}_i\) be the set of strongly connected components of the subgraph induced by \(T\mtp{i}\).
    Let \(\mathcal{N} = \sum_{i \in \positiveintegers} \mathcal{N}_i\).
    Then, 
    \begin{equation*}
		T = \sum_{i \in \positiveintegers} T\mtp{i} = \sum_{i \in \positiveintegers} \sum_{I \in \mathcal{N}_i} I = \sum_{I \in \mathcal{N}} I.
	\end{equation*}

    First, we show that \(\mathcal{N}\) satisfies \ref{item:nestedcollection:disjoint}.
    Let \(I, J \in \mathcal{N}\).
    Then, \(I \in \mathcal{N}_i\) and \(J \in \mathcal{N}_j\) for some \(i, j \in \positiveintegers\).
    Recall that \(I\) is a strongly connected component of the subgraph induced by \(T\mtp{i}\), and \(J\) is a strongly connected component of the subgraph induced by \(T\mtp{j}\).
    Without loss of generality, \(i \leq j\), therefore, \(T\mtp{i} \supset T\mtp{j}\).
    Therefore, there exists a strongly connected component \(K\) of the subgraph induced by \(T\mtp{i}\) such that \(J \subset K\).
    Moreover, since \(I\) and \(K\) are strongly connected components of the subgraph induced by \(T\mtp{i}\), \(I = K\) or \(I \cap K = \varnothing\).
    In the first case, \(I \subset J\), 
    while in the second case, \(I \cap J = \varnothing\),
    as desired.

    Now we show that \(\mathcal{N}\) satisfies \ref{item:nestedcollection:stronglyconnected},
    that is, for any \(\mathcal{R} \subset \mathcal{N}\) such that \(I \cap J = \varnothing\) for all distinct \(I, J \in \mathcal{R}\),
    each \(I \in \mathcal{R}\) is a strongly connected component of the subgraph induced by \(R = \bigcup_{J \in \mathcal{R}} J\).
    We prove by induction on the number of elements of \(\mathcal{R}\).
    If \(\mathcal{R} = \varnothing\), then the statement is vacuously true.
    Suppose that \(\mathcal{R}\) has at least one element.
    Let \(i \in \positiveintegers\) be the minimum index such that \(\mathcal{R} \cap \mathcal{N}_i \neq \varnothing\), and let \(I \in \mathcal{R} \cap \mathcal{N}_i\).
    Then, all \(J \in \mathcal{R}\) are subsets of \(T\mtp{i}\), and consequently, \(R = \bigcup_{J \in \mathcal{R}} J\) is a subset of \(T\mtp{i}\).
    Therefore, \(I \subset R \subset T\mtp{i}\).
    Since \(I\) is a strongly connected component of the subgraph induced by \(T\mtp{i}\), it follows that \(I\) is a strongly connected component of the subgraph induced by \(R\).
    The connected components of the subgraph induced by \(R\) different from \(I\) are the connected components of \(R \setminus I\).
    Applying the induction hypothesis to \(\mathcal{R} \setminus \{I\}\), we obtain that each \(J \in \mathcal{R} \setminus \{I\}\) is a strongly connected component of the subgraph induced by \(R \setminus I\), and therefore, each \(J \in \mathcal{R}\) is a strongly connected component of the subgraph induced by \(R\), as desired.

    Therefore, \(\mathcal{N}\) is a nested multicollection and \(T = \sum_{I \in \mathcal{N}} I\), as desired.
\end{proof}

\begin{proof}[Proof (injectivity)]
    Let \(T\) be a multiset with elements from \(\vertices\). We prove that there is a unique nested multicollection \(\mathcal{N}\) such that \(T = \sum_{I \in \mathcal{N}} I\).
    The proof is by induction on the number of elements of \(T\).
    If \(T = \varnothing\), then \(\mathcal{N} = \varnothing\) is the unique nested multicollection such that \(T = \sum_{I \in \mathcal{N}} I\).

    Suppose that \(T\) has at least one element.
    Let \(\mathcal{N}\) be a nested multicollection such that \(T = \sum_{I \in \mathcal{N}} I\).
    Let \(T^{(1)}\) be the set of elements of \(T\) with multiplicity at least \(1\).
    Let \(\mathcal{R}\) be the set of maximal elements of \(\mathcal{N}\).
    Note that, for all \(v \in T^{(1)}\), there exists a unique \(I \in \mathcal{R}\) such that \(v \in I\).
    Therefore, \(T^{(1)} = \sum_{I \in \mathcal{R}} I\).
    Since \(\mathcal{N}\) is nested and \(I \cap J = \varnothing\) for all distinct \(I, J \in \mathcal{R}\), it follows that each \(I \in \mathcal{R}\) is a strongly connected component of the subgraph induced by \(T^{(1)}\).
    Hence, \(\mathcal{R}\) is the set of strongly connected components of the subgraph induced by \(T^{(1)}\).

    Let \(T' = T - T^{(1)}\), and let \(\mathcal{N}' = \mathcal{N} - \mathcal{R}\).
    Since \(T = \sum_{I \in \mathcal{N}} I\) and \(T^{(1)} = \sum_{I \in \mathcal{R}} I\),
    it follows that \(T' = \sum_{I \in \mathcal{N}'} I\).
    Since \(\mathcal{N}\) is nested, \(\mathcal{N}'\) is nested.
    By the induction hypothesis, \(\mathcal{N}'\) is unique, and consequently, \(\mathcal{N}\) is unique.
\end{proof}

For our discussion of bidirected trees in Section~\ref{section:trees},
we naturally consider nested collections of bidirected graphs.
In this context, condition~\ref{item:nestedcollection:stronglyconnected} of a nested collection can be simplified,
as stated in Lemma~\ref{lemma:nestedcollection:bidirected}.
This lemma is used in the proof of Proposition~\ref{proposition:tree}.

\begin{lemm} \label{lemma:nestedcollection:bidirected}
    Assume \(\graph\) is a bidirected graph.
    A multicollection \(\mathcal{N}\) is a nested multicollection if and only if
    \begin{enumerate}[label = \textup{(N\arabic*')}, ref = \textup{(N\arabic*')}]
        \item[\textup{(N1)}]
        for every pair \(I, J \in \mathcal{N}\),
        either \(I \subset J\), \(J \subset I\), or \(I \cap J = \varnothing\),
        \setcounter{enumi}{1}
        \item \label{item:nestedcollection:stronglyconnected:bidirected}
        for every pair \(I, J \in \mathcal{N}\) such that \(I \cap J = \varnothing\),
        there are no edges between \(I\) and \(J\), and
        \item \label{item:nestedcollection:connected}
        every \(I \in \mathcal{N}\) is connected.
    \end{enumerate}
\end{lemm}

\begin{proof}
    Note that, when \(\graph\) is a bidirected graph,
    strongly connectedness is equivalent to connectedness.
    Let \(\mathcal{N}\) be a multicollection satisfying \ref{item:nestedcollection:disjoint}.
    We show that condition \ref{item:nestedcollection:stronglyconnected} is equivalent to conditions \ref{item:nestedcollection:stronglyconnected:bidirected} and \ref{item:nestedcollection:connected}.

    First, assume that \ref{item:nestedcollection:stronglyconnected} holds.
    Then, applying \ref{item:nestedcollection:stronglyconnected} to \(\mathcal{R} = \{I\}\) for each \(I \in \mathcal{N}\), we obtain \ref{item:nestedcollection:connected}.
    Moreover, applying \ref{item:nestedcollection:stronglyconnected} to \(\mathcal{R} = \{I, J\}\) for each pair \(I, J \in \mathcal{N}\) such that \(I \cap J = \varnothing\), we obtain \ref{item:nestedcollection:stronglyconnected:bidirected}.

    Second, assume that \ref{item:nestedcollection:stronglyconnected:bidirected} and \ref{item:nestedcollection:connected} hold.
    Let \(\mathcal{R} \subset \mathcal{N}\) be such that \(I \cap J = \varnothing\) for all distinct \(I, J \in \mathcal{R}\).
    Then, \ref{item:nestedcollection:stronglyconnected:bidirected} implies that there are no edges between any pair of distinct elements of \(\mathcal{R}\), and \ref{item:nestedcollection:connected} implies that each \(I \in \mathcal{R}\) is connected.
    Therefore, the set of connected components of the subgraph induced by \(\bigcup_{J \in \mathcal{R}} J\) is \(\mathcal{R}\),
    and \ref{item:nestedcollection:stronglyconnected} holds.
\end{proof}

\subsection{Laurent polynomial ring} \label{subsection:laurentpolynomialring}

Let the \vocab{coefficient ring} \(\coefficientring\) be a ring over \(\integers\)
containing elements \(A_v\) for each \(v \in \vertices\) which are algebraically independent.
For example, \(\coefficientring\) could be
\(\integers[A_v : v \in \vertices]\).
Let \(\laurentpolynomialring\) denote the \vocab{Laurent polynomial ring} over \(\coefficientring\)
in the independent variables \(X_v\) for \(v \in \vertices\), that is,
\( \laurentpolynomialring = \coefficientring[X_v^{\pm 1} : v \in \vertices]\).

The monomials in \(\laurentpolynomialring\) in the variables \(X_v\) and \(X_v^{-1}\) for \(v \in \vertices\) are called \vocab{Laurent monomials}.
Any Laurent monomial can be written as
\begin{equation*}
    \laurentmonomial{U}{T} = \prod_{v \in U} X_v \bigg/ \prod_{v \in T} X_v,
\end{equation*}
where \(U\) and \(T\) are disjoint multisets with elements in \(\vertices\).
As a module over \(\coefficientring\), the Laurent polynomial ring \(\laurentpolynomialring\) has a basis consisting of all Laurent monomials.

\subsection{Multifunctions} \label{subsection:multifunctions}

Let \(I\) be a multiset with elements from \(\vertices\).
A \vocab{multifunction} \(f\) (of \(\graph\)) on \(I\) is a directed multigraph with vertex set \(\vertices\) and edge multiset \(\edges_f\) such that
\begin{enumerate}[label = \textup{(F\arabic*)}, ref = \textup{(F\arabic*)}]
    \item for each vertex \(v \in \vertices\), the outdegree of \(v\) in \(f\) is its multiplicity in \(I\), and
    \item each edge in \(\edges_f\) is either a loop or an edge in \(\edges\).
\end{enumerate}
The notation \(f:I\) means that \(f\) is a multifunction on \(I\).
If \(I\) is a set, then a multifunction on \(I\) is naturally identified with a function from \(I\) to \(\vertices\), explaining the chosen nomenclature. 

A multifunction \(f\) is \vocab{acyclic} if the only cycles in \(f\) are loops.
Given two multifunctions \(f\) and \(g\) on multisets \(I\) and \(J\), respectively, we define their \vocab{sum} \(f + g\) as the multifunction on \(I + J\) obtained by taking the sum of the edge multisets of \(f\) and \(g\).

For example, take \(\graph\) to be the directed graph in Figure~\ref{fig:directedgraphs:1}.
Examples of multifunctions of \(\graph\) are given in Figure~\ref{fig:multifunctions}.
The multifunctions in Figure~\ref{fig:multifunctions:1} and Figure~\ref{fig:multifunctions:2} are acyclic, and the multifunction in Figure~\ref{fig:multifunctions:3} is not acyclic.

\begin{figure}[htpb]
    \centering
    \begin{subfigure}[t]{0.29\textwidth}
        \centering
        \begin{tikzpicture}
            \node[shape=circle, draw=black] (2) at (0,0){2};
            \node[shape=circle, draw=black] (3) at (-1,-1){3};
            \node[shape=circle, draw=black] (1) at (-1,1){1};
            \node[shape=circle, draw=black] (4) at (-2,0){4};
            \path[->]
                (1) edge[bend left] node {} (2)
                (1) edge[bend right] node {} (2)
                (2) edge[] node {} (3)
                (1) edge[loop above] node {} (1);
        \end{tikzpicture}.
        \caption{A multifunction on \(\{1,\allowbreak 1,\allowbreak 1,\allowbreak 2\}\) with edge multiset \(\{(1, 1),\allowbreak (1, 2),\allowbreak (1, 2),\allowbreak (2, 3)\}\).}
        \label{fig:multifunctions:1}
    \end{subfigure}
    \hfill
    \begin{subfigure}[t]{0.28\textwidth}
        \centering
        \begin{tikzpicture}
            \node[shape=circle, draw=black] (2) at (0,0){2};
            \node[shape=circle, draw=black] (3) at (-1,-1){3};
            \node[shape=circle, draw=black] (1) at (-1,1){1};
            \node[shape=circle, draw=black] (4) at (-2,0){4};
            \path[->]
                (3) edge[] node {} (4)
                (4) edge[] node {} (1);
        \end{tikzpicture}.
        \caption{A multifunction on \(\{3,\allowbreak 4\}\) with edge multiset \(\{(3, 4),\allowbreak (4, 1)\}\).}
        \label{fig:multifunctions:2}
    \end{subfigure}
    \hfill
    \begin{subfigure}[t]{0.29\textwidth}
        \centering
        \begin{tikzpicture}
            \node[shape=circle, draw=black] (2) at (0,0){2};
            \node[shape=circle, draw=black] (3) at (-1,-1){3};
            \node[shape=circle, draw=black] (1) at (-1,1){1};
            \node[shape=circle, draw=black] (4) at (-2,0){4};
            \path[->]
                (1) edge[bend left] node {} (2)
                (1) edge[bend right] node {} (2)
                (2) edge[] node {} (3)
                (1) edge[loop above] node {} (1)
                (3) edge[] node {} (4)
                (4) edge[] node {} (1);
        \end{tikzpicture}.
        \caption{The sum of the multifunctions in Figures~\ref{fig:multifunctions:1} and~\ref{fig:multifunctions:2}.}
        \label{fig:multifunctions:3}
    \end{subfigure}
    \caption{Multifunctions of the directed graph in Figure~\ref{fig:directedgraphs:1}.}
    \label{fig:multifunctions}
\end{figure}

Given a multifunction \(f\) and a subset \(W \subset V\),
the \vocab{restriction} of \(f\) to \(W\), denoted by \(f|_{W}\), is the multifunction whose edge multiset consists of the edges in \(\edges_f\) with source in \(W\).

\subsection{Weight of a multifunction} \label{subsection:weightofamultifunction}

The \vocab{weight} of a multifunction \(f\) on \(I\), denoted by \(\weight(f)\), is the non-Laurent polynomial in \(\laurentpolynomialring\) given by
\begin{equation*}
    \weight(f)
    = \prod_{(v, w) \in \edges_f} \Tilde{X}_{(v, w)},
\end{equation*}
where
\begin{equation*}
    \Tilde{X}_{(v, w)} = 
    \begin{cases}
    X_{w} & \text{if } w \neq v, \\
    A_{v} & \text{if } w = v.
    \end{cases}
\end{equation*}
For example,
the weights of the multifunctions in Figures~\ref{fig:multifunctions:1}, \ref{fig:multifunctions:2}, and \ref{fig:multifunctions:3} are, respectively,
\begin{equation*}
    A_1 X_2^2 X_3,
    \qquad
    X_1 X_4,
    \quad
    \text{and}
    \qquad
    A_1 X_1 X_2^2 X_3 X_4.
\end{equation*}

The \vocab{normalized weight} of a multifunction \(f\) on \(I\), denoted by \(\normalizedweight(f)\), is the Laurent polynomial in \(\laurentpolynomialring\) given by

\begin{equation*}
    \normalizedweight(f)
    = \frac{\weight(f)}{\prod_{v \in I} X_v}
    = \prod_{(v, w) \in \edges_f} \frac{\Tilde{X}_{(v, w)}}{X_v}.
\end{equation*}
For example,
the normalized weights of the multifunctions in Figures~\ref{fig:multifunctions:1}, \ref{fig:multifunctions:2}, and \ref{fig:multifunctions:3} are, respectively,
\begin{equation*}
    \frac{A_1 X_2^2 X_3}{X_1^3 X_2} = 
    \frac{A_1 X_2 X_3}{X_1^3},
    \qquad
    \frac{X_1 X_4}{X_3 X_4}=
    \frac{X_1}{X_3},
    \quad
    \text{and}
    \qquad
    \frac{A_1 X_1 X_2^2 X_3 X_4}{X_1^3 X_2 X_3 X_4} =
    \frac{A_1 X_2}{X_1^2}.
\end{equation*}

We remark that the nomenclature of ``weights'' and ``normalized weights'' is not used by \cite{linearLP}, although the concept is present in their work.

Note that both weights and normalized weights are products over edges of the multifunction, and therefore, if \(f\) and \(g\) are multifunctions, then 
\begin{equation*}
    \weight(f + g) = \weight(f) \weight(g)
    \qquad \text{and} \qquad
    \normalizedweight(f + g) = \normalizedweight(f) \normalizedweight(g).
\end{equation*}
Note that the normalized weight of a cycle is \(1\).
This is an important observation that is used in Lemma~\ref{lemma:set:decomposition:acyclic} and Proposition~\ref{proposition:alongpath:normalizedweight}.

\subsection{The \texorpdfstring{\(Y\)}{Y} Laurent polynomials} \label{subsection:Ylaurentpolynomials}

Let \(I\) be a subset of \(\vertices\).
We define the Laurent polynomial \(Y_I \in \laurentpolynomialring\) by
\begin{equation*}
    Y_I = \frac{\sumacy{f : I} \weight(f)}{\prod_{i \in I} X_i},
\end{equation*}
where the sum is over all acyclic multifunctions \(f\) on \(I\).
We write \(Y_I\) in terms of normalized weights as \(Y_I = \sumacy{f:I} \normalizedweight(f)\).

Recall that, since \(I\) is a set, a multifunction on \(I\) is naturally identified with a function from \(I\) to \(\vertices\).

For example,
taking \(\graph\) to be the directed graph in Figure~\ref{fig:directedgraphs:1},
the Laurent polynomial \(Y_{\{1, 2\}}\) is
\begin{equation*}
    \medmuskip=2mu
    \frac{
        A_1X_1 + X_3X_1 + X_4X_1 +
        A_1A_2 + X_2A_2 + X_3A_2 + X_4A_2 +
        A_1X_3 + X_2X_3 + X_3^2  + X_4X_3
    }{X_1X_2}.
\end{equation*}
Note that the eleven terms in the numerator correspond to the eleven acyclic multifunctions on \(\{1, 2\}\),
which are all twelve multifunctions on \(\{1, 2\}\) by assigning one of the possible four edges to \(1\),
and one of the possible three edges to \(2\),
except for the assignment of the edge \((1, 2)\) to \(1\) and the edge \((2, 1)\) to \(2\) which is not acyclic.

\subsection{Graph LP algebra and clusters} \label{subsection:graphLPalgebraandclusters}

The central algebraic structure of our research is the \vocab{graph Laurent phenomenon algebra} \(\graphLPalgebra{\graph}\) over \(\coefficientring\) associated to \(\graph\), defined by \cite{linearLP}.
Lam and Pylyavskyy~\cite{linearLP} proved that \(\graphLPalgebra{\graph}\) is the algebra over \(\coefficientring\) generated by 
\(X_v\) for \(v \in \vertices\) and \(Y_I\) for strongly connected \(I \subset \vertices\).
The reader who is new to graph LP algebras can take this as the definition of \(\graphLPalgebra{\graph}\).

Although only the Laurent polynomials \(Y_I\) for strongly connected \(I \subset \vertices\) are generators of \(\graphLPalgebra{\graph}\),
the Laurent polynomials \(Y_I\) for \(I \subset \vertices\) are in \(\graphLPalgebra{\graph}\) as well, as guaranteed by Lemma~\ref{lemma:y:product}.

\begin{lemm}[{\cite[Lemma 4.2]{linearLP}}] \label{lemma:y:product}
    Let \(I\) be a subset of \(\vertices\).
    Then, \(Y_I = \prod_{J} Y_J\), where the product is over the strongly connected components \(J\) of the subgraph of \(\graph\) induced by \(I\).
\end{lemm}

The monomials in \(\graphLPalgebra{\graph}\) in the elements \(X_v\) for \(v \in v\) and \(Y_I\) for \(I \subset \vertices\) are simply called \vocab{monomials}.
Any monomial in \(\graphLPalgebra{\graph}\) can be written as
\begin{equation*}
    \monomial{U}{\mathcal{S}} =
    \prod_{v \in U} X_v \prod_{I \in \mathcal{S}} Y_I,
\end{equation*}
where \(U\) is a multiset with elements in \(\vertices\) and \(\mathcal{S}\) is a multiset with elements in \(\powerset{\vertices}\).
The set of monomials in \(\graphLPalgebra{\graph}\) is a spanning set of \(\graphLPalgebra{\graph}\) as a module over \(\coefficientring\).
The monomials in \(\graphLPalgebra{\graph}\) in the elements \(Y_I\) for \(I \subset \vertices\) are called \vocab{\(Y\)-monomials}, which can be written as
\( \monomial{\varnothing}{\mathcal{S}} =
    \prod_{I \in \mathcal{S}} Y_I
\)
where \(\mathcal{S}\) is a multiset with elements in \(\powerset{\vertices}\).

The elements \(X_v\) for \(v \in \vertices\) and \(Y_I\) for strongly connected \(I \subset \vertices\) are grouped into sets called \vocab{clusters}.
We refer to \cite{linearLP} for the definition of clusters.
Lam and Pylyavskyy~\cite{linearLP} proved that the clusters of \(\graphLPalgebra{\graph}\) are the sets of the form
\begin{equation*}
    \{ X_v : v \in U\} \cup \{Y_I : I \in \mathcal{N}\}
\end{equation*}
where \(U \subseteq \vertices\) and \(\mathcal{N}\) is a maximal nested collection on \(\vertices \setminus U\).
The reader who is new to graph LP algebras can take this as the definition of the clusters of \(\graphLPalgebra{\graph}\).

The monomials in \(\graphLPalgebra{\graph}\) in the elements of a given cluster are called \vocab{cluster monomials}.
Any cluster monomial can be written as
\begin{equation*}
    \monomial{U}{\mathcal{N}} =
    \prod_{v \in U} X_v \prod_{I \in \mathcal{N}} Y_I,
\end{equation*}
where \(U\) is a multiset with elements in \(\vertices\) and \(\mathcal{N}\) is a nested collection on \(\vertices \setminus U\mtp{1}\).
The monomials in \(\graphLPalgebra{\graph}\) in the elements \(Y_I\) for \(I \subset \vertices\) of a given cluster are called \vocab{cluster \(Y\)-monomials}, which can be written as \(\monomial{\varnothing}{\mathcal{N}} = \prod_{I \in \mathcal{N}} Y_I\) where \(\mathcal{N}\) is a nested collection on \(\vertices\).

Recall that Lemma~\ref{lemma:multinested:unique} provides a bijection between nested multicollections \(\mathcal{N}\) and multisets \(T\) with elements from \(\vertices\).
Therefore, we can reindex the cluster monomials by a pair \((U, T)\) where \(U\) and \(T\) are disjoint multisets with vertices from \(\vertices\),
and the cluster monomial indexed by \((U, T)\) is 
\( \monomial{U}{T} = \monomial{U}{\mathcal{N}} \),
where \(\mathcal{N}\) is the unique nested multicollection such that \(T = \sum_{I \in \mathcal{N}} I\).

\section{Linear independence} \label{section:linearindependence}

In this section, we show that the cluster monomials in a graph LP algebra are linearly independent, as stated in Theorem~\ref{theorem:linearindependence}.

\begin{theo} \label{theorem:linearindependence}
    The set of cluster monomials of a graph LP algebra is linearly independent over \(\coefficientring\).
\end{theo}

The overview of the proof of Theorem~\ref{theorem:linearindependence} is as follows.
The key observation is that the cluster monomial \(\monomial{U}{T}\) is a linear combination of the Laurent monomials \( \laurentmonomial{U'}{T'} \) with \(U' \subseteq U\) and \(T' \supseteq T\).
This observation induces a triangular structure on the coefficient matrix relating cluster monomials to Laurent monomials, with the diagonal entries being nonzero.
From this triangular structure, the linear independence of the cluster monomials follows from the linear independence of the Laurent monomials.

Lemmas~\ref{lemma:coefficient:U1T1U2T2} and~\ref{lemma:coefficient:UT} establish the triangularity condition,
and then we prove Theorem~\ref{theorem:linearindependence}.
Recall that \(\laurentmonomial{U}{T}\) denotes the Laurent monomial \({\prod_{v \in U} X_v} \big/ {\prod_{v \in T}X_v}\). 

\begin{lemm} \label{lemma:coefficient:U1T1U2T2}
    Let \(U_1, U_2, T_1, T_2\) be multisets with vertices from \(\vertices\)
    such that \(U_1 \cap T_1 = U_2 \cap T_2 = \varnothing\).
    If \(U_1 \not\subset U_2\), then, when expanding \(\monomial{U_1}{T_1}\) as a linear combination of Laurent monomials over \(\coefficientring\), the coefficient of \(\laurentmonomial{U_2}{T_2}\) is \(0\).
    If \(T_1 \not\supset T_2\), then, when expanding \(\monomial{U_1}{T_1}\) as a linear combination of Laurent monomials over \(\coefficientring\), the coefficient of \(\laurentmonomial{U_2}{T_2}\) is \(0\).
\end{lemm}

\begin{proof}
    Let \(\mathcal{N}_1\) be the unique nested multicollection such that \(T_1 = \sum_{I \in \mathcal{N}_1} I\).
    Recall that
    \begin{equation*}
        \monomial{U_1}{T_1}
        = \prod_{v \in U_1} X_v \prod_{I \in \mathcal{N}_1} Y_I
        = \frac{\prod_{v \in U_1} X_v}{\prod_{v \in T_1} X_v}
        \left(
            \prod_{I \in \mathcal{N}_1} \sumacy{f \colon I} \weight(f)
        \right).
    \end{equation*} 

    Since \( \prod_{I \in \mathcal{N}_1} \sumacy{f \colon I} \weight(f) \) is a polynomial, the numerator of any term of \(\monomial{U_1}{T_1}\) is a multiple of \(\prod_{v \in U_1} X_v\).
    Hence, if \(U_1 \not\subset U_2\),
    then the coefficient of
    \begin{equation*}
		\laurentmonomial{U_2}{T_2} = \frac{\prod_{v \in U_2} X_v}{\prod_{v \in T_2} X_v}
	\end{equation*}
    in \(\monomial{U_1}{T_1}\) is \(0\).

    Since \( \prod_{I \in \mathcal{N}_1} \sumacy{f \colon I} \weight(f) \) is a polynomial,
    the denominator of any term of \(\monomial{U_1}{T_1}\) is a divisor of \(\prod_{v \in T_1} X_v\).
    Hence, if \(T_1 \not\supset T_2\),
    then the coefficient of
    \begin{equation*}
		\laurentmonomial{U_2}{T_2} = \frac{\prod_{v \in U_2} X_v}{\prod_{v \in T_2} X_v}
	\end{equation*}
    in \(\monomial{U_1}{T_1}\) is \(0\).
\end{proof}

\begin{lemm} \label{lemma:coefficient:UT}
    Let \(U, T\) be multisets with vertices from \(\vertices\).
    When expanding \(\monomial{U}{T}\) as a linear combination of Laurent monomials over \(\coefficientring\), the coefficient of \(\laurentmonomial{U}{T}\) is \(\prod_{v \in T} A_v\).
\end{lemm}

\begin{proof}
    Let \(\mathcal{N}\) be the unique nested multicollection such that \(T = \sum_{I \in \mathcal{N}} I\).
    Recall that
    \begin{equation*}
        \monomial{U}{T}
        = \prod_{v \in U} X_v \prod_{I \in \mathcal{N}} Y_I 
        = \frac{\prod_{v \in U} X_v}{\prod_{v \in T} X_v}
        \left(
            \prod_{I \in \mathcal{N}} \sumacy{f \colon I} \weight(f)
        \right).
    \end{equation*}

    Note that \( \prod_{I \in \mathcal{N}} \sumacy{f \colon I} \weight(f) \) is a polynomial with constant term \(\prod_{v \in T} A_v\), obtained from the constant acyclic function on each \(I \in \mathcal{N}\).
    Therefore, the coefficient of \(\laurentmonomial{U}{T}\) in \(\monomial{U}{T} = {\prod_{v \in U} X_v} \big/ {\prod_{v \in T}X_v}\) is \(\prod_{v \in T} A_v\).
\end{proof}

We are now ready to prove Theorem~\ref{theorem:linearindependence}.

\begin{proof}[Proof of Theorem~\ref{theorem:linearindependence}]
    Suppose that there exists a nontrivial linear combination of cluster monomials that equals \(0\).
    Explicitly, suppose that there exists coefficients \(\coefficient{U}{T}\), not all zero, such that
    \begin{equation*}
        \sum_{U, T} \coefficient{U}{T} \monomial{U}{T} = 0.
    \end{equation*}

    Pick \((U_0, T_0)\) with \(\coefficient{U_0}{T_0} \neq 0\) such that \(|U_0| - |T_0|\) is minimized.
    This implies that, for all \((U, T) \neq (U_0, T_0)\) such that \(\coefficient{U}{T} \neq 0\),
    it holds that \(|U| - |T| \geq |U_0| - |T_0|\), and consequently, \(U \not\subset U_0\) or \(T \not\supset T_0\).

    Hence, by Lemma~\ref{lemma:coefficient:U1T1U2T2}, the coefficient of \(\laurentmonomial{U_0}{T_0}\) in 
    \( \coefficient{U}{T} \monomial{U}{T} \)
    is \(0\)
    for all \((U, T) \neq (U_0, T_0)\) such that \(\coefficient{U}{T} \neq 0\);
    and, by Lemma~\ref{lemma:coefficient:UT}, the coefficient of \(\laurentmonomial{U_0}{T_0}\) in the expression
    \( \coefficient{U_0}{T_0} \monomial{U_0}{T_0} \)
    is 
    \( \coefficient{U_0}{T_0} \prod_{v \in T_0} A_v \).
    Therefore, the coefficient of \(\laurentmonomial{U_0}{T_0}\) in
    \( \sum_{U, T} \coefficient{U}{T} \monomial{U}{T} \)
    is
    \( \coefficient{U_0}{T_0} \prod_{v \in T_0} A_v \neq 0\),
    a contradiction.
\end{proof}

We remark that the triangularity argument used to prove linear independence does not directly imply that the cluster monomials form a spanning set.
In finite-dimensional vector spaces, any order of the indices of the basis has a minimal element,
and the triangularity argument would imply that the linear combinations of the basis elements form a spanning set.
However, this intuition does not fully carry over to infinite-dimensional vector spaces.
In our context, the issue arises from the absence of a minimal element in the ordering of pairs \((U, T)\) of multisets with vertices from \(\vertices\).

\section{Spanning set} \label{section:spanningset}

In this section, we show that the cluster monomials in a graph LP algebra form a spanning set, as stated in Theorem~\ref{theorem:spanningset}.

\begin{theo} \label{theorem:spanningset}
    Any monomial is a linear combination of cluster monomials over \(\coefficientring\).
\end{theo}

The overview of the proof of Theorem~\ref{theorem:spanningset} is as follows.
We show that any \(Y\)-monomial can be expressed as a linear combination of cluster \(Y\)-monomials over \(\integers\).
This involves establishing a relation on the normalized weights of tuples of functions on multicollections with equal sum.
Using the fact that the normalized weight of a cycle is \(1\),
we derive a corresponding relation on acyclic functions,
which we then translate into a relation on \(Y\)-monomials.
We then use \cite[Lemma~4.7]{linearLP}, a result on the product of an \(X\)-variable and a \(Y\)-variable, to show Theorem~\ref{theorem:spanningset}.

\subsection{Product of an \texorpdfstring{\(X\)}{X}-variable and a \texorpdfstring{\(Y\)}{Y}-variable} \label{subsection:productXY}

Let \(p \colon \ppath{v}{I}{w}\) denote the statement that \(p\) is a vertex non-repeating directed path from \(v\) to \(w\) with intermediary vertices in \(I\).
Let \(I \setminus p\) denote the set of vertices in \(I\) that are not in \(p\).
Lemma~\ref{lemma:XYgood} allows us to rewrite the product of \(X_v\) and \(Y_I\) whenever \(v \in I\).

\begin{lemm}[{\cite[inferred from Lemma 4.7]{linearLP}}] \label{lemma:XYgood}
Let \(I \in \powerset{\vertices}\) and \(v \in I\).
Then,
\begin{equation*} \label{ppathsequation}
    X_v Y_I =
    \sum_{w \in I}
    \sum_{p}^{\ppath{v}{I}{w}}
    Y_{I \setminus p}  X_w
    +
    \sum_{w \not\in I}
    \sum_{p}^{\ppath{v}{I}{w}}
    Y_{I \setminus p}  A_w.
\end{equation*}
\end{lemm}

Any monomial only in the \(X\)-variables is already a cluster monomial, and Lemma \ref{lemma:XYgood} essentially gives us a way to multiply an \(X\)-variable with a \(Y\)-variable.
It remains to understand how to multiply \(Y\)-variables. 

Thus, our proof strategy is to focus first on \(Y\)-monomials. We show first that cluster \(Y\)-monomials span \(Y\)-monomials. Then, we use Lemma \ref{lemma:XYgood} to deduce that cluster monomials span all monomials. 

\subsection{Integer expansion of \texorpdfstring{\(Y\)}{Y}-monomials} \label{subsection:Ymonomialsspanning}

In this subsection, we show that cluster \(Y\)-monomials span \(Y\)-monomials, as stated in Theorem~\ref{theorem:Yclusters-span}.

\begin{theo} \label{theorem:Yclusters-span}
    Any \(Y\)-monomial is an integer linear combination of cluster \(Y\)-monomials.
\end{theo}

Before proving Theorem~\ref{theorem:Yclusters-span} in its generality,
we discuss a special case that motivates the proof.
Consider the \(Y\)-monomial \(\monomial{\varnothing}{\{I, J\}} = Y_I Y_J\)
where \(I\) and \(J\) are \emph{disjoint} subsets.
Note that 
\begin{equation*}
    Y_I Y_J = \sum_{f, g} \normalizedweight(f + g),
\end{equation*}
where the sum is over acyclic functions \(f\) on \(I\) and \(g\) on \(J\).
Given such a pair \(f, g\), their sum \(f + g\) is a function on \(I \cup J\).
However, \(f + g\) may not be acyclic,
as cycles that are not entirely contained in \(I\) or in \(J\) may appear in \(f + g\).
Since the normalized weight of a cycle is \(1\),
removing these cycles do not affect the normalized weight of \(f + g\).
With this in mind, one might expect that
\begin{equation*}
    Y_I Y_J = \sum_{C} Y_{I \cup J \setminus C},
\end{equation*}
where the sum runs over families \(C\) of disjoint cycles in \(I \cup J\)
that are not fully contained in \(I\) or \(J\) (including the empty family),
and \(V(C)\) denotes the set of vertices in the cycles in \(C\).
Finally, the right-hand-side is the sum of a cluster \(Y\)-monomial \(Y_{I \cup J}\)
and \emph{lower order terms} \(Y_{I \cup J \setminus C}\) that can be dealt with inductively.

The more general setting has two complications:
the monomial \(\monomial{\varnothing}{\mathcal{S}}\) might have \(|\mathcal{S}| > 2\);
and the sets in \(\mathcal{S}\) have nontrivial intersections.
Nevertheless, the core idea is roughly the same:
study the sum of acyclic functions on the sets in \(\mathcal{S}\),
how they relate to functions on the sum of the sets in \(\mathcal{S}\),
and use induction to handle lower order terms.

Given a tuple of multifunctions, we define its normalized weight as the product of the normalized weights of its elements.

\begin{lemm}[Preimages Lemma] \label{lemma:preimage}
    Let \(\mathcal{F}_1, \mathcal{F}_2, \mathcal{G}\) be sets of tuples of multifunctions.
    If \(\phi_1 \colon \mathcal{F}_1 \to \mathcal{G}\) and
    \(\phi_2 \colon \mathcal{F}_2 \to \mathcal{G}\) preserve normalized weights, and
    \( |\phi_1^{-1}(g)| = |\phi_2^{-1}(g)| \) for all \(g \in \mathcal{G}\), then
    \begin{equation*}
        \sum_{f \in \mathcal{F}_1} \normalizedweight(f) = \sum_{f \in \mathcal{F}_2} \normalizedweight(f).
    \end{equation*} 
\end{lemm}

\begin{proof}
    For \(i \in \{1, 2\}\), apply the weight-preserving property and double-count pairs \((f, g) \in \mathcal{F} \times \mathcal{G}\) such that \(\phi_i(f) = g\) to obtain
    \begin{equation*}
        \sum_{f \in \mathcal{F}} \normalizedweight(f) =
        \sum_{f \in \mathcal{F}} \normalizedweight(\phi_i(f)) =
        \sum_{g \in \mathcal{G}} \left| \phi_i^{-1}(g) \right| \normalizedweight(g).
    \end{equation*} 
    Therefore,
    \begin{equation*}
        \sum_{f \in \mathcal{F}_1} \normalizedweight(f) =
        \sum_{g \in \mathcal{G}} \left| \phi_1^{-1}(g) \right| \normalizedweight(g) = 
        \sum_{g \in \mathcal{G}} \left| \phi_2^{-1}(g) \right| \normalizedweight(g) = 
        \sum_{f \in \mathcal{F}_2} \normalizedweight(f). \qedhere
    \end{equation*} 
\end{proof}

\begin{lemm} \label{lemma:multiset:decomposition}
    Let \(\mathcal{S}\) be a multiset with elements in \(\powerset{\vertices}\).
    Let \(T = \sum_{I \in \mathcal{S}} I\) be the multiset containing the vertices that appear in the sets of \(\mathcal{S}\), counting multiplicities.
    Let \(\mathcal{T}\) be the multiset whose elements are \(T\mtp{1}, T\mtp{2}, \dots \in \powerset{\vertices}\).
    Then,
    \begin{equation*} \label{eq:multiset:decomposition:nwt}
        \prod_{S \in \mathcal{S}} \sum_{f : S} \normalizedweight(f) = \prod_{T' \in \mathcal{T}} \sum_{f : T'} \normalizedweight(f).
    \end{equation*}
\end{lemm}

\begin{proof}
    Let \(S_1, S_2, \dots, S_{|\mathcal{S}|}\) be the elements of \(\mathcal{S}\).
    Let \(T_1, T_2, \dots, T_{|\mathcal{T}|}\) be the elements of \(\mathcal{T}\).
    Note that \(T = \sum_{i=1}^{|\mathcal{S}|} S_i = \sum_{j=1}^{|\mathcal{T}|} T_j\).
    Define the sets
    \begin{gather*}
        \mathcal{F}_1 = \{ (f_1, f_2, \dots, f_{|\mathcal{S}|}) : f_i \text{ is a multifunction on } S_i \}, \\
        \mathcal{F}_2 = \{ (f_1, f_2, \dots, f_{|\mathcal{T}|}) : f_j \text{ is a multifunction on } T_j \}, \\
        \mathcal{G} = \{ g : g \text{ is a multifunction on } T \},
    \end{gather*}
    and the \(\normalizedweight\)-preserving functions
    \( \phi_1 \colon \mathcal{F}_1 \to \mathcal{G} \) and
    \( \phi_2 \colon \mathcal{F}_2 \to \mathcal{G} \)
    defined by \(\phi_1(f_1, f_2, \allowbreak \dots, f_{|\mathcal{S}|}) = f_1 + f_2 + \cdots + f_{|\mathcal{S}|}\) and \(\phi_2(f_1, f_2, \dots, f_{|\mathcal{T}|}) = f_1 + f_2 + \cdots + f_{|\mathcal{T}|}\).

    Fix \(g \in \mathcal{G}\).
    Let's count the number of \((f_1, f_2, \dots, f_{|\mathcal{S}|}) \in \mathcal{F}_1\) such that \(\phi_1(f_1, f_2, \allowbreak \dots, f_{|\mathcal{S}|}) = g\).
    Fix \(v \in \vertices\).
    Let \(g_v\) denote the multiset of edges in \(g\) from \(v\).
    Note that \(|g_v|\) is the multiplicity of \(v\) in \(T\),
    which is the number of sets \(S_i\) in \(\mathcal{S}\) that contain \(v\).
    Hence, the number of ways to appropriately assign the edges of \(g_v\) to the functions \(f_1\), \(f_2\), \(\dots\), \(f_{|\mathcal{S}|}\) is
    \begin{equation*}
        \frac{|g_v|!}{1^{\left| g_v\mtp{1} \right|} \cdot 2^{\left| g_v\mtp{2} \right|} \cdot 3^{\left| g_v\mtp{3} \right|} \cdots}.
    \end{equation*}
    Therefore, the number of \((f_1, f_2, \dots, f_{|\mathcal{S}|}) \in \mathcal{F}_1\) such that \(\phi_1(f_1, f_2, \dots, f_{|\mathcal{S}|}) = g\) is
    \begin{equation*}
        \left| \phi_1^{-1}(g) \right| =
        \prod_{v \in \vertices} \frac{\left| g_v \right|!}{1^{\left| g_v\mtp{1} \right|} \cdot 2^{\left| g_v\mtp{2} \right|} \cdot 3^{\left| g_v\mtp{3} \right|} \cdots}.
    \end{equation*}

    The same argument applies to \(\phi_2\),
    and we obtain, for all \(g \in \mathcal{G}\),
    \begin{equation*}
        \left| \phi_1^{-1}(g) \right| =
        \prod_{v \in \vertices} \frac{\left| g_v \right|!}{1^{\left| g_v\mtp{1} \right|} \cdot 2^{\left| g_v\mtp{2} \right|} \cdot 3^{\left| g_v\mtp{3} \right|} \cdots} =
        \left| \phi_2^{-1}(g) \right|.
    \end{equation*}
    Hence, by Lemma~\ref{lemma:preimage}, we obtain 
    \begin{equation*}
        \sum_{f \in \mathcal{F}_1} \normalizedweight(f)
        = \sum_{f \in \mathcal{F}_2} \normalizedweight(f). \qedhere
    \end{equation*}
\end{proof}

Given \(S \in \mathcal{P}(V)\),
we denote by \(\mathcal{C}_S\) the set of families of vertex-disjoint simple directed cycles in the restriction of \(\graph\) to \(S\).
Given a family \(C \in \mathcal{C}_S\),
we denote by \(V(C)\) the set of vertices that appear in the cycles in \(C\).
For example, if \(S = \{1, 2, 3, 4\}\) and \(\graph\) is as in Figure~\ref{fig:directedgraphs:1},
then \(\mathcal{C}_S\) contains \(14\) elements:
\(1\) empty family,
\(4\) families of a single \(3\)-cycle,
\(2\) families of a single \(4\)-cycle,
\(5\) families of a single \(2\)-cycle,
and \(2\) families of two \(2\)-cycles.
These \(14\) families are shown in Figure~\ref{fig:familyofcycles}.

\begin{figure}[htbp]
    \centering
    \begin{tikzpicture}[scale=0.6]
        \node[shape=circle, draw=black, inner sep=1pt] (2) at (0,0){2};
        \node[shape=circle, draw=black, inner sep=1pt] (3) at (-1,-1){3};
        \node[shape=circle, draw=black, inner sep=1pt] (1) at (-1,1){1};
        \node[shape=circle, draw=black, inner sep=1pt] (4) at (-2,0){4};
    \end{tikzpicture}
    \hfill
    \begin{tikzpicture}[scale=0.6]
        \node[shape=circle, draw=black, inner sep=1pt] (2) at (0,0){2};
        \node[shape=circle, draw=black, inner sep=1pt] (3) at (-1,-1){3};
        \node[shape=circle, draw=black, inner sep=1pt] (1) at (-1,1){1};
        \node[shape=circle, draw=black, inner sep=1pt] (4) at (-2,0){4};
        \path[draw,thick, ->]
            (1) edge node {} (2)
            (2) edge node {} (3)
            (3) edge node {} (1);
    \end{tikzpicture}
    \hfill
    \begin{tikzpicture}[scale=0.6]
        \node[shape=circle, draw=black, inner sep=1pt] (2) at (0,0){2};
        \node[shape=circle, draw=black, inner sep=1pt] (3) at (-1,-1){3};
        \node[shape=circle, draw=black, inner sep=1pt] (1) at (-1,1){1};
        \node[shape=circle, draw=black, inner sep=1pt] (4) at (-2,0){4};
        \path[draw,thick, ->]
            (1) edge node {} (3)
            (3) edge node {} (2)
            (2) edge node {} (1);
    \end{tikzpicture}
    \hfill
    \begin{tikzpicture}[scale=0.6]
        \node[shape=circle, draw=black, inner sep=1pt] (2) at (0,0){2};
        \node[shape=circle, draw=black, inner sep=1pt] (3) at (-1,-1){3};
        \node[shape=circle, draw=black, inner sep=1pt] (1) at (-1,1){1};
        \node[shape=circle, draw=black, inner sep=1pt] (4) at (-2,0){4};
        \path[draw,thick, ->]
            (1) edge node {} (3)
            (3) edge node {} (4)
            (4) edge node {} (1);
    \end{tikzpicture}
    \hfill
    \begin{tikzpicture}[scale=0.6]
        \node[shape=circle, draw=black, inner sep=1pt] (2) at (0,0){2};
        \node[shape=circle, draw=black, inner sep=1pt] (3) at (-1,-1){3};
        \node[shape=circle, draw=black, inner sep=1pt] (1) at (-1,1){1};
        \node[shape=circle, draw=black, inner sep=1pt] (4) at (-2,0){4};
        \path[draw,thick, ->]
            (1) edge node {} (4)
            (4) edge node {} (3)
            (3) edge node {} (1);
    \end{tikzpicture}
    \hfill
    \begin{tikzpicture}[scale=0.6]
        \node[shape=circle, draw=black, inner sep=1pt] (2) at (0,0){2};
        \node[shape=circle, draw=black, inner sep=1pt] (3) at (-1,-1){3};
        \node[shape=circle, draw=black, inner sep=1pt] (1) at (-1,1){1};
        \node[shape=circle, draw=black, inner sep=1pt] (4) at (-2,0){4};
        \path[draw,thick, ->] (1) edge node {} (2);
        \path[draw,thick, ->] (2) edge node {} (3);
        \path[draw,thick, ->] (3) edge node {} (4);
        \path[draw,thick, ->] (4) edge node {} (1);
    \end{tikzpicture}
    \hfill
    \begin{tikzpicture}[scale=0.6]
        \node[shape=circle, draw=black, inner sep=1pt] (2) at (0,0){2};
        \node[shape=circle, draw=black, inner sep=1pt] (3) at (-1,-1){3};
        \node[shape=circle, draw=black, inner sep=1pt] (1) at (-1,1){1};
        \node[shape=circle, draw=black, inner sep=1pt] (4) at (-2,0){4};
        \path[draw,thick, ->] (1) edge node {} (4);
        \path[draw,thick, ->] (2) edge node {} (1);
        \path[draw,thick, ->] (3) edge node {} (2);
        \path[draw,thick, ->] (4) edge node {} (3);
    \end{tikzpicture}
    \\[.5em]
    \begin{tikzpicture}[scale=0.6]
        \node[shape=circle, draw=black, inner sep=1pt] (2) at (0,0){2};
        \node[shape=circle, draw=black, inner sep=1pt] (3) at (-1,-1){3};
        \node[shape=circle, draw=black, inner sep=1pt] (1) at (-1,1){1};
        \node[shape=circle, draw=black, inner sep=1pt] (4) at (-2,0){4};
        \path[draw,thick, ->, bend left] (1) edge node {} (2);
        \path[draw,thick, ->, bend left] (2) edge node {} (1);
    \end{tikzpicture}
    \hfill
    \begin{tikzpicture}[scale=0.6]
        \node[shape=circle, draw=black, inner sep=1pt] (2) at (0,0){2};
        \node[shape=circle, draw=black, inner sep=1pt] (3) at (-1,-1){3};
        \node[shape=circle, draw=black, inner sep=1pt] (1) at (-1,1){1};
        \node[shape=circle, draw=black, inner sep=1pt] (4) at (-2,0){4};
        \path[draw,thick, ->, bend left] (1) edge node {} (3);
        \path[draw,thick, ->, bend left] (3) edge node {} (1);
    \end{tikzpicture}
    \hfill
    \begin{tikzpicture}[scale=0.6]
        \node[shape=circle, draw=black, inner sep=1pt] (2) at (0,0){2};
        \node[shape=circle, draw=black, inner sep=1pt] (3) at (-1,-1){3};
        \node[shape=circle, draw=black, inner sep=1pt] (1) at (-1,1){1};
        \node[shape=circle, draw=black, inner sep=1pt] (4) at (-2,0){4};
        \path[draw,thick, ->, bend left] (1) edge node {} (4);
        \path[draw,thick, ->, bend left] (4) edge node {} (1);
    \end{tikzpicture}
    \hfill
    \begin{tikzpicture}[scale=0.6]
        \node[shape=circle, draw=black, inner sep=1pt] (2) at (0,0){2};
        \node[shape=circle, draw=black, inner sep=1pt] (3) at (-1,-1){3};
        \node[shape=circle, draw=black, inner sep=1pt] (1) at (-1,1){1};
        \node[shape=circle, draw=black, inner sep=1pt] (4) at (-2,0){4};
        \path[draw,thick, ->, bend left] (2) edge node {} (3);
        \path[draw,thick, ->, bend left] (3) edge node {} (2);
    \end{tikzpicture}
    \hfill
    \begin{tikzpicture}[scale=0.6]
        \node[shape=circle, draw=black, inner sep=1pt] (2) at (0,0){2};
        \node[shape=circle, draw=black, inner sep=1pt] (3) at (-1,-1){3};
        \node[shape=circle, draw=black, inner sep=1pt] (1) at (-1,1){1};
        \node[shape=circle, draw=black, inner sep=1pt] (4) at (-2,0){4};
        \path[draw,thick, ->, bend left] (4) edge node {} (3);
        \path[draw,thick, ->, bend left] (3) edge node {} (4);
    \end{tikzpicture}
    \hfill
    \begin{tikzpicture}[scale=0.6]
        \node[shape=circle, draw=black, inner sep=1pt] (2) at (0,0){2};
        \node[shape=circle, draw=black, inner sep=1pt] (3) at (-1,-1){3};
        \node[shape=circle, draw=black, inner sep=1pt] (1) at (-1,1){1};
        \node[shape=circle, draw=black, inner sep=1pt] (4) at (-2,0){4};
        \path[draw,thick, ->, bend left] (1) edge node {} (2);
        \path[draw,thick, ->, bend left] (2) edge node {} (1);
        \path[draw,thick, ->, bend left] (3) edge node {} (4);
        \path[draw,thick, ->, bend left] (4) edge node {} (3);
    \end{tikzpicture}
    \hfill
    \begin{tikzpicture}[scale=0.6]
        \node[shape=circle, draw=black, inner sep=1pt] (2) at (0,0){2};
        \node[shape=circle, draw=black, inner sep=1pt] (3) at (-1,-1){3};
        \node[shape=circle, draw=black, inner sep=1pt] (1) at (-1,1){1};
        \node[shape=circle, draw=black, inner sep=1pt] (4) at (-2,0){4};
        \path[draw,thick, ->, bend left] (1) edge node {} (4);
        \path[draw,thick, ->, bend left] (4) edge node {} (1);
        \path[draw,thick, ->, bend left] (2) edge node {} (3);
        \path[draw,thick, ->, bend left] (3) edge node {} (2);
    \end{tikzpicture}
    \caption{The families of vertex-disjoint simple directed cycles in \(\mathcal{C}_S\) for \(S = \{1, 2, 3, 4\}\) and \(\graph\) as in Figure~\ref{fig:directedgraphs:1}.
    In order, the set \(V(C)\) of vertices in each family \(C\) is
    \(\emptyset\),
    \(\{1, 2, 3\}\), \(\{1, 2, 3\}\), \(\{1, 3, 4\}\), \(\{1, 3, 4\}\),
    \(\{1, 2, 3, 4\}\), \(\{1, 2, 3, 4\}\),
    \(\{1, 2\}\), \(\{1, 3\}\), \(\{1, 4\}\), \(\{2, 3\}\), \(\{3, 4\}\),
    \(\{1, 2, 3, 4\}\), \(\{1, 2, 3, 4\}\).}
    \label{fig:familyofcycles}
\end{figure}

\begin{lemm} \label{lemma:set:decomposition:acyclic}
    Let \(S \in \powerset{\vertices}\).
    Then,
    \begin{equation*}
        \sum_{f : S} \normalizedweight(f)
        = \sum_{C \in \mathcal{C}_S} Y_{S \setminus V(C)}.
    \end{equation*}
\end{lemm}

\begin{proof}
    Given a multifunction \(f\) on \(S\),
    since \(S\) is a set,
    the outdegree of a vertex in \(f\) is at most \(1\) and,
    consequently, no two cycles in \(f\) share a vertex.
    Therefore,
    each multifunction \(f \colon S\) can be uniquely decomposed as the sum of a family \(C\) of vertex-disjoint cycles and an acyclic function \(g \colon S \setminus V(C)\).
    Similarly, given a family of vertex-disjoint cycles \(C\) and an acyclic function \(g \colon S \setminus V(C)\),
    a function \(f \colon S\) can be obtained by taking the sum of \(g\) and the cycles in \(C\).
    Thus, there is a bijection between the set of functions \(f \colon S\) and
    the set of pairs \((C, g)\) where \(C\) is a family of vertex-disjoint cycles in \(S\) and \(g\) is an acyclic function on \(S \setminus V(C)\).
    Moreover,
    if this bijection associates \(f\) to \((C, g)\),
    then the normalized weight of \(f\) is the product of the normalized weight of \(g\) and the normalized weight of the cycles in \(C\).
    Since the normalized weight of a cycle is \(1\), it follows that \(\normalizedweight(f) = \normalizedweight(g)\).
    Finally, we have
    \begin{equation*}
        \sum_{f : S} \normalizedweight(f)
        = \sum_{C \in \mathcal{C}_S} \sumacy{g : S \setminus V(C)} \normalizedweight(g)
        = \sum_{C \in \mathcal{C}_S} Y_{S \setminus V(C)}. \qedhere
    \end{equation*}
\end{proof}

\begin{lemm} \label{lemma:multiset:cluster}
    Let \(T\) be a multiset of vertices.
    Let \(\mathcal{T}\) be the multiset whose elements are \(T\mtp{1}, T\mtp{2}, \dots \in \powerset{\vertices}\).
    Then, the \(Y\)-monomial \(\monomial{\varnothing}{\mathcal{T}}\) is a cluster \(Y\)-monomial.
\end{lemm}

\begin{proof}
    Let \(\mathcal{N}_i\) be the set of strongly connected components of the subgraph of \(\graph\) induced by \(T\mtp{i}\) for each \(i \in \positiveintegers\),
    and let \(\mathcal{N} = \sum_i \mathcal{N}_i\).
    Lemma~\ref{lemma:multinested:unique} implies that \(\mathcal{N}\) is a nested multicollection.
    Finally, Lemma~\ref{lemma:y:product} implies that \(Y_{T\mtp{i}} = \prod_{J \in \mathcal{N}_i} Y_J\),
    and consequently
    \begin{equation*}
        \monomial{\varnothing}{\mathcal{T}}
        = \prod_{i} Y_{T\mtp{i}}
        = \prod_{i} \prod_{J \in \mathcal{N}_i} Y_J
        = \prod_{J \in \mathcal{N}} Y_J
        = \monomial{\varnothing}{\mathcal{N}}. \qedhere
    \end{equation*}
\end{proof}

With these lemmas in hand, we can prove the main result of this section, Theorem~\ref{theorem:Yclusters-span}.

\begin{theo*}[Theorem~\ref{theorem:Yclusters-span}, repeated]
    Any \(Y\)-monomial is an integer linear combination of cluster \(Y\)-monomials.

    More precisely, for any multiset \(\mathcal{S}\) with elements in \(\powerset{\vertices}\), the \(Y\)-monomial \(\monomial{\varnothing}{\mathcal{S}}\) is an integer linear combination of cluster \(Y\)-monomials.
    Moreover, a cluster \(Y\)-monomials \(\monomial{\varnothing}{\mathcal{N}}\) appearing in the integer linear combination of \(\monomial{\varnothing}{\mathcal{S}}\) satisfies that all vertices in a set of \(\mathcal{N}\) are vertices in a set of \(\mathcal{S}\).
\end{theo*}

\begin{proof}
    We apply strong induction on \(\left| \sum_{I \in \mathcal{S}} I \right|\),
    that is, the number of vertices in a set of \(\mathcal{S}\), counting multiplicities.
    Assume, by induction hypothesis, that the result holds for any multiset \(\mathcal{R}\) with elements in \(\powerset{\vertices}\) such that \(\left| \sum_{J \in \mathcal{R}} J \right| < \left| \sum_{I \in \mathcal{S}} I \right|\).

    Let \(T = \sum_{I \in \mathcal{S}} I\), the multiset containing the vertices in a set of \(\mathcal{S}\), counting multiplicities.
    Let \(\mathcal{T}\) be the multiset whose elements are \(T\mtp{1}, T\mtp{2}, \dots \in \powerset{\vertices}\).
    By Lemma~\ref{lemma:multiset:cluster}, the \(Y\)-monomial \(\monomial{\varnothing}{\mathcal{T}}\) is a cluster \(Y\)-monomial.

    If \(\mathcal{T} = \mathcal{S}\), then \(\monomial{\varnothing}{\mathcal{S}} = \monomial{\varnothing}{\mathcal{T}}\) is a cluster \(Y\)-monomial, and we are done.
    Otherwise, \(\mathcal{T} \neq \mathcal{S}\).
    From Lemma~\ref{lemma:multiset:decomposition}, we know that
    \begin{equation*}
        \prod_{S \in \mathcal{S}} \sum_{f : S} \normalizedweight(f) = \prod_{T \in \mathcal{T}} \sum_{f : T} \normalizedweight(f).
    \end{equation*}
    Applying Lemma~\ref{lemma:set:decomposition:acyclic} to each term of both products, we obtain
    \begin{equation*}
        \prod_{S \in \mathcal{S}} \sum_{C_S \in \mathcal{C}_S} Y_{S \setminus C_S}
        =
        \prod_{T \in \mathcal{T}} \sum_{C_{T} \in \mathcal{C}_{T}} Y_{T \setminus C_{T}}.
    \end{equation*}
    Rewriting the equation above in terms of \(Y\)-monomials, we obtain
    \begin{equation*}
        \sum_{ \substack{
                (C_S\mtp{1}\!\!,\,
                C_S\mtp{2}\!\!,\,\dots) \\
                \rotatebox{-90}{\(\in\)} \\
                \mathcal{C}_S\mtp{1}\!\times \mathcal{C}_S\mtp{2}\!\times \cdots }
                }
                \monomial{\varnothing}{
                    \{ 
                        S\mtp{1} \setminus C_S\mtp{1},\,
                        \dots
                    \}
                }
        =
        \sum_{ \substack{
                (C_{T}\mtp{1}\!\!,\,
                C_{T}\mtp{2}\!\!,\,\dots) \\
                \rotatebox{-90}{\(\in\)} \\
                \mathcal{C}_{T}\mtp{1}\!\times \mathcal{C}_{T}\mtp{2}\!\times \cdots }
                }
                \monomial{\varnothing}{
                    \{ 
                        T\mtp{1} \setminus C_{T}\mtp{1},\,
                        \dots
                    \}
                }.
    \end{equation*}
    Note that for \(\varnothing = C_S\mtp{1} = C_S\mtp{2} = \cdots\), the term in the sum in the left-hand side is the \(Y\)-monomial \(\monomial{\varnothing}{\mathcal{S}}\).
    Therefore, the \(Y\)-monomial \(\monomial{\varnothing}{\mathcal{S}}\) is equal to
    \begin{equation} \label{eq:rewrite-mos}
        \sum_{ \substack{
            (C_{T}\mtp{1}\!\!,\,
            C_{T}\mtp{2}\!\!,\,\dots) \\
            \rotatebox{-90}{\(\in\)} \\
            \mathcal{C}_{T}\mtp{1}\!\times \mathcal{C}_{T}\mtp{2}\!\times \cdots }
            }
            \mkern -9mu
            \monomial{\varnothing}{
                \{ 
                    T\mtp{1} \setminus C_{T}\mtp{1},\,
                    \dots
                \}
            }
        -
        \sum_{ \substack{
            (C_S\mtp{1}\!\!,\,
            C_S\mtp{2}\!\!,\,\dots) \\
            \rotatebox{-90}{\(\in\)} \\
            \mathcal{C}_S\mtp{1}\!\times \mathcal{C}_S\mtp{2}\!\times \cdots \\
            \text{except all }\varnothing} 
            }
            \mkern -9mu
            \monomial{\varnothing}{
                \{ 
                    S\mtp{1} \setminus C_S\mtp{1},\,
                    \dots
                \}
            }.
    \end{equation}

    Note that each term in \eqref{eq:rewrite-mos} is of the form \(\pm \monomial{\varnothing}{\mathcal{R}}\),
    and satisfies either that \(\mathcal{R} = \mathcal{T}\) or that \(\left|\sum_{I \in \mathcal{R}} I\right| < \left|\sum_{I \in \mathcal{S}} I\right|\).
    Hence, by Lemma~\ref{lemma:multiset:cluster} and by the induction hypothesis, each monomial \(\monomial{\varnothing}{\mathcal{R}}\) in \eqref{eq:rewrite-mos} is an integer linear combination of cluster \(Y\)-monomials, and therefore, the \(Y\)-monomial \(\monomial{\varnothing}{\mathcal{S}}\) is an integer combination of cluster \(Y\)-monomials.

    Moreover, note that each monomial in \eqref{eq:rewrite-mos} of the form \(\pm \monomial{\varnothing}{\mathcal{R}}\) have the property that all vertices in a set of \(\mathcal{R}\) are vertices in a set of \(\mathcal{S}\).
    Consequently, by the induction hypothesis, the cluster \(Y\)-monomials \(\monomial{\varnothing}{\mathcal{N}}\) in the expansion of \(\monomial{\varnothing}{\mathcal{S}}\) are such that all vertices in a set of \(\mathcal{N}\) are vertices in a set of \(\mathcal{S}\).

    Therefore, by induction, the result holds.
\end{proof}

\subsection{Expansion of monomials} \label{subsection:proofofspanningset}

We proceed to the proof of the main result of this section, Theorem~\ref{theorem:spanningset}.

\begin{theo*}[Theorem~\ref{theorem:spanningset}, repeated] \label{theorem:spanningset:repeated}
    Any monomial is a linear combination of cluster monomials over \(\coefficientring\).

    More precisely, for any multiset \(U\) with elements in \(\vertices\) and any multiset \(\mathcal{S}\) with elements in \(\powerset{\vertices}\), the monomial \(\monomial{U}{\mathcal{S}}\) is a linear combination of cluster monomials over \(\coefficientring\).
\end{theo*}

\begin{proof}
    We apply strong induction on \(\left|\sum_{I \in \mathcal{S}} I\right|\), the number of vertices in a set of \(\mathcal{S}\), counting multiplicities.
    Assume, by induction hypothesis, that the monomial \(\monomial{T}{\mathcal{R}}\) is a linear combination over \(\coefficientring\) of cluster monomials whenever \(\left|\sum_{J \in \mathcal{R}} J\right| < \left|\sum_{I \in \mathcal{S}} I\right|\).

    Assume there exists \(v \in U\) such that \(v \in S\) for some \(S \in \mathcal{S}\).
    We can write
    \begin{equation*}
        \monomial{U}{\mathcal{S}} = \monomial{U - \{v\}}{\mathcal{S} - \{S\}} \cdot X_v \cdot Y_S.
    \end{equation*}
    Apply Lemma~\ref{lemma:XYgood} to \(X_v\) and \(Y_S\) to obtain
    \begin{equation*}
        X_v Y_S =
        \sum_{w \in S}
        \sum_{p}^{\ppath{v}{S}{w}}
        Y_{S \setminus p}  X_w
        +
        \sum_{w \not\in S}
        \sum_{p}^{\ppath{v}{S}{w}}
        Y_{S \setminus p}  A_w.
    \end{equation*}
    Therefore,
    \begin{align*}
        \monomial{U}{\mathcal{S}}
        = & \sum_{w \in S}
        \sum_{p}^{\ppath{v}{S}{w}}
        \monomial{U - \{v\} + \{w\}}{\mathcal{S} - \{S\} + \{S \setminus p\}} \\
        & + 
        \sum_{w \not\in S}
        \sum_{p}^{\ppath{v}{S}{w}}
        A_w \monomial{U - \{v\}}{\mathcal{S} + \{S \setminus p\} - \{S\}}.
    \end{align*}

    This implies that we can write \(\monomial{U}{\mathcal{S}}\) as a linear combination over \(\coefficientring\) of monomials of the form \(\monomial{T}{\mathcal{R}}\),
    where \(\left|\sum_{J \in \mathcal{R}} J\right| < \left|\sum_{I \in \mathcal{S}} I\right|\).
    Therefore, by the induction hypothesis, \(\monomial{U}{\mathcal{S}}\) is a linear combination over \(\coefficientring\) of cluster monomials.

    Otherwise, assume all \(v \in U\) satisfy \(v \not\in S\) for all \(S \in \mathcal{S}\).
    In other words, \(\mathcal{S} \in \powerset{\vertices \setminus U}\).
    We can write
    \( \monomial{U}{\mathcal{S}} = \monomial{U}{\varnothing} \monomial{\varnothing}{\mathcal{S}} \).
    Apply Theorem~\ref{theorem:Yclusters-span} to obtain that \( \monomial{\varnothing}{\mathcal{S}} \) is a linear combination over \(\coefficientring\) of cluster \(Y\)-monomials of the form \(\monomial{\varnothing}{\mathcal{N}}\), where \(\mathcal{N}\) is a nested multicollection such that the vertices in sets of \(\mathcal{N}\) are a subset of the vertices in sets of \(\mathcal{S}\).
    Explicitly,
    \begin{equation*}
        \monomial{\varnothing}{\mathcal{S}} = \sum_{\mathcal{N}} c(\mathcal{N}) \monomial{\varnothing}{\mathcal{N}},
    \end{equation*}
    where the sum is over all nested multicollections \(\mathcal{N}\) such that the vertices in sets of \(\mathcal{N}\) are a subset of the vertices in sets of \(\mathcal{S}\).
    Multiplying by \(\monomial{U}{\varnothing}\), we obtain
    \begin{equation*}
        \monomial{U}{\mathcal{S}} = \sum_{\mathcal{N}} c(\mathcal{N}) \monomial{U}{\mathcal{N}}.
    \end{equation*}
    Since \(\mathcal{S} \in \powerset{\vertices \setminus U}\), it follows that \(\mathcal{N} \in \powerset{\vertices \setminus U}\) for all \(\mathcal{N}\) in the sum above, and consequently, \(\monomial{U}{\mathcal{N}}\) is a cluster monomial.
    Therefore, \(\monomial{U}{\mathcal{S}}\) is a linear combination over \(\coefficientring\) of cluster monomials.
\end{proof}

\section{Nonnegativity for trees}\label{section:trees}

In this section, we prove Theorem~\ref{theorem:tree}, which is the special case of Conjecture~\ref{conjecture:positivity} where \(\graph\) is a bidirected tree.
Throughout this section, we assume that \(\graph\) is a bidirected tree,
that is, a directed graph such that \(vu \in E(\graph) \iff uv \in E(\graph)\) and whose underlying undirected graph is a tree.

\begin{theo*}[Theorem~\ref{theorem:tree}, repeated]
    If \(\graph\) is a bidirected tree, each monomial of \(\graphLPalgebra{\graph}\) is a linear combination of cluster monomials over \(\coefficientring\) with nonnegative coefficients.
\end{theo*}

The outline of the proof of Theorem~\ref{theorem:tree} is as follows.
First, in Proposition~\ref{proposition:expansionfortrees},
we explicitly expand \(Y_IY_J\) as a sum of \(Y_{I \cup J}Y_{I \cap J}\) and \emph{lower order terms} of the form \(Y_{K}Y_{L}\) with \(K \subseteq I \cup J\) and \(L \subseteq I \cap J\).
This is done similarly to the proof of Theorem~\ref{theorem:Yclusters-span},
but the restriction to bidirected trees and to the product of two \(Y\)-variables allows us to derive an explicit formula that guarantees nonnegative coefficients.
Then, using that \(\graph\) is bidirected, it follows that the \(Y\)-monomials of \(\graphLPalgebra{\graph}\) are nonnegative integer linear combinations of cluster \(Y\)-monomials.
Finally, Theorem~\ref{theorem:tree} follows using Lemma~\ref{lemma:XYgood}.

We remark that in order to prove Conjecture \ref{conjecture:positivity} where \(\graph\) is a bidirected graph that is not necessarily a tree, the only missing step is to generalize Proposition~\ref{proposition:expansionfortrees} by showing that, if \(\graph\) is a bidirected graph, then \(Y_IY_J\) is a sum of \(Y_{I \cup J}Y_{I \cap J}\) and lower order terms.

\subsection{Nonnegative formula for the product of two \texorpdfstring{\(Y\)}{Y}-variables}

Lam and Pylyavskyy \cite[Theorem 6.1]{linearLP} give a Ptolemy-like formula for expanding \(Y_IY_J\) into a positive integer linear combination of cluster monomials in the case where \(\graph\) is a path.
In order to extend their formula, we introduce the notion of a path from \(I\) to \(J\).

\begin{defi}[Path from \(I\) to \(J\)]
    Let \(\graph\) be a bidirected tree.
    Let \(I, J \subset \vertices\).
    A path \(P = w_1w_2\dots w_k\) of \(\graph\) is said to be a \vocab{path from \(I\) to \(J\)}
    if \(w_1 \in I \setminus J\), \(w_2, \dots, w_{k-1} \in I \cap J\), and \(w_k \in J \setminus I\).

    Let \(W_I(P)\), \(W_J(P)\), \(W_{I \cup J}(P)\), and \(W_{I \cap J}(P)\)
    be the sets of vertices of \(P\) that are in \(I\), in \(J\), in \(I \cup J\), and in \(I \cap J\), respectively.
    Explicitly,
    \begin{align*}
        W_I(P) & = \{w_1, w_2, w_3, \dots, w_{k-1}\}, &
        W_J(P) & = \{w_2, w_3, \dots, w_k\}, \\
        W_{I \cup J}(P) & = \{w_1, w_2, w_3, \dots, w_k\}, &
        W_{I \cap J}(P) & = \{w_2, w_3, \dots, w_{k-1}\}.
    \end{align*}

    Given a family \(\mathcal{P} = \{P_1, P_2, \dots, P_t\}\) of disjoint paths from \(I\) to \(J\),
    define \(W_X(\mathcal{P}) = \bigcup_{i=1}^t W_X(P_i)\) for \(X \in \{I, J, I \cup J, I \cap J\}\),
    that is, the set of vertices of the paths in \(\mathcal{P}\) that are in \(X\) for \(X \in \{I, J, I \cup J, I \cap J\}\).
\end{defi}

Note that the empty set is a family of (zero) disjoint paths from \(I\) to \(J\).

\begin{figure}[htbp]
    \centering
    \begin{tikzpicture}[x=1cm]
        \draw[fill=blue!40, draw=black, fill opacity=.5] (-.5, -.7) rectangle (5.5, 1.6);
        \draw[fill=red!40, draw=black, fill opacity=.5] (.5, -.6) rectangle (6.5, 1.7);
        \node[shape=circle, draw=black] (1) at (0,0){};
        \node[shape=circle, draw=black] (2) at (1,0){};
        \node[shape=circle, draw=black] (3) at (2,0){};
        \node[] (4) at (3,0){\small \(\cdots\)};
        \node[shape=circle, draw=black] (5) at (4,0){};
        \node[shape=circle, draw=black] (6) at (5,0){};
        \node[shape=circle, draw=black] (7) at (6,0){};
        \node[shape=circle, draw=black] (1b) at (0,1){};
        \node[shape=circle, draw=black] (2b) at (1,1){};
        \node[shape=circle, draw=black] (3b) at (2,1){};
        \node[] (4b) at (3,1){\small \(\cdots\)};
        \node[shape=circle, draw=black] (5b) at (4,1){};
        \node[shape=circle, draw=black] (6b) at (5,1){};
        \node[shape=circle, draw=black] (7b) at (6,1){};
        \draw[very thick] (1) -- (2) -- (3) -- (4) -- (5) -- (6) -- (7);
        \draw[very thick] (1b) -- (2b) -- (3b) -- (4b) -- (5b) -- (6b) -- (7b);
    \end{tikzpicture}
    \caption{Illustration of a family \(\mathcal{P}\) of (two) disjoint paths from \(I\) to \(J\).
    The vertices of the blue region are in \(W_I(\mathcal{P})\),
    the vertices of the red region are in \(W_J(\mathcal{P})\),
    the vertices of blue or red regions are in \(W_{I \cup J}(\mathcal{P})\),
    and the vertices of blue and red regions are in \(W_{I \cap J}(\mathcal{P})\).}
    \label{fig:disjointpaths}
\end{figure}

Proposition~\ref{proposition:expansionfortrees} gives an explicit formula for expanding \(Y_IY_J\) into a positive integer linear combination of cluster monomials.

\begin{prop} \label{proposition:expansionfortrees}
    Let \(\graph\) be a tree.
    Let \(I, J \subset \vertices\).
    The expansion of \(Y_IY_J\) into cluster monomials is given by
    \begin{equation*}
        Y_IY_J = \sum_{\mathcal{P}} Y_{I \cup J \setminus W_{I \cup J}(\mathcal{P})}Y_{I \cap J \setminus W_{I \cap J}(\mathcal{P})},
    \end{equation*} 
    where the sum is taken over all families \(\mathcal{P}\) of disjoint paths from \(I\) to \(J\).
\end{prop}

Before proving Proposition \ref{proposition:expansionfortrees},
we give two examples.

\begin{exam}
    Let \(\graph = P_n\), the path graph with vertices \(\{1, 2, \dots, n\}\).
    Figure~\ref{fig:tree:1} shows the path graph \(P_6\).
    Let \(I = \{1, 2, \dots, k\}\) and \(J = \{l, l+1, \dots, n\}\).
    If \(k < l - 1\), then there is no path from \(I\) to \(J\) and Proposition~\ref{proposition:expansionfortrees} implies that \(Y_{I} Y_{J} = Y_{I \cup J}\), which is a cluster monomial.
    If \(k \geq l - 1\), then the only path from \(I\) to \(J\) is \(P=\{l-1, l, \dots, k, k+1\}\).
    Therefore, \(Y_{I} Y_{J}\) is expanded into two cluster monomials:
    one corresponding to the empty family of paths,
    and the other corresponding to the family containing only \(P\).
    Explicitly, we recover \cite[Theorem~6.1]{linearLP}'s Ptolemy-like formula:
    \begin{equation*}
        Y_{\{1, \dots, k\}} Y_{\{l, \dots, n\}}
        = Y_{\{1, \dots, n\}} Y_{\{l, \dots, k\}}
          + Y_{\{1, \dots, l-2, k+2, \dots, n\}}.
    \end{equation*}
\end{exam}

\begin{figure}[htbp]
    \centering
    \begin{subfigure}[b]{0.47\textwidth}
        \centering
        \begin{tikzpicture}
            \useasboundingbox (-.5,-.6) rectangle (5.5,.6);
            \node[shape=circle, draw=black] (1) at (0,0){1};
            \node[shape=circle, draw=black] (2) at (1,0){2};
            \node[shape=circle, draw=black] (3) at (2,0){3};
            \node[shape=circle, draw=black] (4) at (3,0){4};
            \node[shape=circle, draw=black] (5) at (4,0){5};
            \node[shape=circle, draw=black] (6) at (5,0){6};
            \path[draw,thick]
                (1) edge node {} (2)
                (2) edge node {} (3)
                (3) edge node {} (4)
                (4) edge node {} (5)
                (5) edge node {} (6);
        \end{tikzpicture}
        \caption{The path graph \(P_6\), with six vertices, \(1\), \(2\), \(3\), \(4\), \(5\), and \(6\), and five edges, \(12\), \(23\), \(34\), \(45\), and \(56\).}
        \label{fig:tree:1}
    \end{subfigure}
    \hfill
    \begin{subfigure}[b]{0.47\textwidth}
        \centering
        \begin{tikzpicture}
            \node[shape=circle, draw=black] (1) at (0,0){1};
            \node[shape=circle, draw=black] (2) at (1,0){2};
            \node[shape=circle, draw=black] (3) at (2,0){3};
            \node[shape=circle, draw=black] (4) at (0,-1){4};
            \node[shape=circle, draw=black] (5) at (1,-1){5};
            \node[shape=circle, draw=black] (6) at (2,-1){6};
            \path[draw,thick]
                (2) edge node {} (3)
                (1) edge node {} (2)
                (4) edge node {} (5)
                (5) edge node {} (6)
                (5) edge node {} (2);
        \end{tikzpicture}
        \caption{The tree with six vertices, \(1\), \(2\), \(3\), \(4\), \(5\), and \(6\), and five edges, \(23\), \(12\), \(45\), \(56\), and \(52\).}
        \label{fig:tree:2}
    \end{subfigure}
    \caption{Two examples of trees.}
    \label{fig:tree}
\end{figure}

\begin{exam}
    Let \(\graph\) be the tree with vertex set \(\vertices = \{1,\allowbreak 2,\allowbreak 3,\allowbreak 4,\allowbreak 5,\allowbreak 6\}\) and undirected edge set \(\edges = \{23,\allowbreak 12,\allowbreak 45,\allowbreak 56,\allowbreak 52\}\), as shown in Figure~\ref{fig:tree:2}.
    Let \(I = \{1, 2, 4, 5\}\) and let \(J = \{2, 3, 5, 6\}\).
    There are six families of disjoint paths from \(I\) to \(J\), 
    \begin{gather*}
        \varnothing, \quad
        \big\{ (1, 2, 3) \big\}, \quad
        \big\{ (4, 5, 6) \big\}, \quad
        \big\{ (1, 2, 5, 6) \big\}, \quad
        \big\{ (4, 5, 2, 3) \big\}, \quad
        \big\{ (1, 2, 3), (4, 5, 6) \big\}.
    \end{gather*}
    Therefore, Proposition~\ref{proposition:expansionfortrees} implies that \(Y_IY_J\) is expanded into six cluster monomials, one for each family of disjoint paths. Explicitly,
    \begin{equation*}
        Y_{\{1, 2, 4, 5\}} Y_{\{2, 3, 5, 6\}}
        = Y_{\{1, 2, 3, 4, 5, 6\}} Y_{\{2, 5\}}
        + Y_{\{4, 5, 6\}} Y_{\{5\}}
        + Y_{\{1, 2, 3\}} Y_{\{2\}}
        + Y_{\{4, 3\}}
        + Y_{\{1, 6\}}
        + 1.
    \end{equation*}
\end{exam}

For organizational purposes, the proof of Proposition~\ref{proposition:expansionfortrees} is divided into multiple lemmas.
Before we proceed with the lemmas,
one more definition is salient.

\begin{defi}[Multifunction along path]
    Let \(P = w_1w_2\dots w_k\) be a path from \(I\) to \(J\).
    The \vocab{multifunction along the path \(P\)} is the multifunction \(p\) on the multiset
    \begin{equation*}
        \{w_1, w_2^2, w_3^2, \dots, w_k\} = W_I(P) \sqcup W_J(P) = W_{I \cup J}(P) \sqcup W_{I \cap J}(P)
    \end{equation*}
    with directed edges
    \begin{equation*}
        w_1w_2, w_2w_3, \dots, w_{k-1}w_k,
        w_kw_{k-1}, w_{k-1}w_{k-2}, \dots, w_2w_1.
    \end{equation*}
    Let \(\mathcal{P} = \{P_1, \dots, P_m\}\) be a family of disjoint paths from \(I\) to \(J\).
    The \vocab{multifunction along the family of disjoint paths \(\mathcal{P}\)} is the sum of the multifunctions along each path in \(\mathcal{P}\),
    that is, \(p = \sum_{i=1}^m p_i\).
\end{defi}

\begin{figure}[htbp]
    \centering
    \begin{tikzpicture}[x=1cm]
        \draw[fill=blue!40, draw=black, fill opacity=.5] (-.5, -.7) rectangle (5.5, 1.6);
        \draw[fill=red!40, draw=black, fill opacity=.5] (.5, -.6) rectangle (6.5, 1.7);
        \node[shape=circle, draw=black] (1) at (0,0){};
        \node[shape=circle, draw=black] (2) at (1,0){};
        \node[shape=circle, draw=black] (3) at (2,0){};
        \node[] (4) at (3,0){\small \(\cdots\)};
        \node[shape=circle, draw=black] (5) at (4,0){};
        \node[shape=circle, draw=black] (6) at (5,0){};
        \node[shape=circle, draw=black] (7) at (6,0){};
        \node[shape=circle, draw=black] (1b) at (0,1){};
        \node[shape=circle, draw=black] (2b) at (1,1){};
        \node[shape=circle, draw=black] (3b) at (2,1){};
        \node[] (4b) at (3,1){\small \(\cdots\)};
        \node[shape=circle, draw=black] (5b) at (4,1){};
        \node[shape=circle, draw=black] (6b) at (5,1){};
        \node[shape=circle, draw=black] (7b) at (6,1){};
        \draw[->] (1) edge[bend left] (2);
        \draw[->] (2) edge[bend left] (1);
        \draw[->] (2) edge[bend left] (3);
        \draw[->] (3) edge[bend left] (2);
        \draw[->] (3) edge[bend left] (4);
        \draw[->] (4) edge[bend left] (3);
        \draw[->] (4) edge[bend left] (5);
        \draw[->] (5) edge[bend left] (4);
        \draw[->] (5) edge[bend left] (6);
        \draw[->] (6) edge[bend left] (5);
        \draw[->] (6) edge[bend left] (7);
        \draw[->] (7) edge[bend left] (6);
        \draw[->] (1b) edge[bend left] (2b);
        \draw[->] (2b) edge[bend left] (1b);
        \draw[->] (2b) edge[bend left] (3b);
        \draw[->] (3b) edge[bend left] (2b);
        \draw[->] (3b) edge[bend left] (4b);
        \draw[->] (4b) edge[bend left] (3b);
        \draw[->] (4b) edge[bend left] (5b);
        \draw[->] (5b) edge[bend left] (4b);
        \draw[->] (5b) edge[bend left] (6b);
        \draw[->] (6b) edge[bend left] (5b);
        \draw[->] (6b) edge[bend left] (7b);
        \draw[->] (7b) edge[bend left] (6b);
    \end{tikzpicture}
    \caption{Illustration of the multifunction along the family of disjoint paths \(\mathcal{P}\) from \(I\) to \(J\) in Figure~\ref{fig:disjointpaths}.}
\end{figure}

Since a multifunction along a path from \(I\) to \(J\) consists of a sum of \(2\)-cycles, whose normalized weight is \(1\),
we have the following proposition:

\begin{prop} \label{proposition:alongpath:normalizedweight}
    The normalized weight of the multifunction along a family of disjoint paths from \(I\) to \(J\) is \(1\).
\end{prop}

Assume \(\graph\) is a directed subgraph of a bidirected tree.
Let \(I, J \subset \vertices\).
Let \(\mathcal{F}\) be the set of pairs \((f, g)\) such that \(f\) is an acyclic function on \(I\) and \(g\) is an acyclic function on \(J\).
Hence, \(Y_I Y_J = \sum_{(f, g) \in \mathcal{F}} \normalizedweight(f + g)\).

The outline of the proof of Proposition \ref{proposition:expansionfortrees} is as follows.
After proving some properties of multifunctions along paths,
we partition \(\mathcal{F}\) into sets \(\mathcal{F}_{\mathcal{P}}\) for each family \(\mathcal{P}\) of disjoint paths from \(I\) to \(J\),
so that the multifunction \(p\) along \(\mathcal{P}\) is a submultiset of the sum \(f + g\).
Then, we show that \(\mathcal{F}_{\mathcal{P}}\) is in a \(\normalizedweight\)-preserving bijection with the set \(\mathcal{H}_{\mathcal{P}}\) of pairs \((f', g')\) such that \(f'\) is an acyclic function on \((I \cup J) \setminus W_{I \cup J}(\mathcal{P})\) and \(g'\) is an acyclic function on \((I \cap J) \setminus W_{I \cap J}(\mathcal{P})\).
Finally, aggregating over all families of disjoint paths, we obtain the desired expansion.

\subsubsection{Properties of multifunctions along paths}

Note that \(I \sqcup J = (I \cup J) + (I \cap J)\).
That is, vertices \(v \in I \cap J\) have multiplicity \(2\) in \(I \sqcup J\),
while vertices \(v \in (I \cup J) \setminus (I \cap J)\) have multiplicity \(1\) in \(I \sqcup J\).

\begin{prop} \label{proposition:restriction-strongly-connected}
    Let \(h\) be a multifunction on a submultiset of \(I \sqcup J\).
    Let \(W \subset I \cup J\) be a strongly connected component of \(h\).
    Assume that \(W\) is not entirely contained in \(I\) or \(J\).
    Then, the restriction \(h|_W\) to edges with source in \(W\)
     is a multifunction along a path from \(I\) to \(J\).
\end{prop}

\begin{proof}
    Let \(k = |W|\).
    Since \(W\) is not entirely contained in \(I\),
    there exists \(w_1 \in W\) such that \(w_1 \in I \setminus J\).
    Since \(W\) is not entirely contained in \(J\),
    there exists \(w_k \in W\) such that \(w_k \in J \setminus I\).

    Consider the undirected graph \(\gamma\) obtained from the subgraph of \(\graph\) induced by \(W\) by removing the directions of the edges.
    Since \(W\) is a strongly connected component of \(h \subset \graph\),
    the undirected graph \(\gamma\) is connected.
    Since \(\graph\) is a subgraph of a bidirected tree,
    the undirected graph \(\gamma\) is a tree.
    Each vertex \(w \in W\) has multiplicity \(2\) in \(I \sqcup J\), except for \(w_1\) and \(w_k\) which have multiplicity \(1\).
    Therefore, each vertex \(w \in W\) of \(\gamma\) has degree at most \(2\), and the vertices \(w_1\) and \(w_k\) have degree \(1\).

    A walk from \(w_1\) to \(w_k\) in \(\gamma\) exists, since \(\gamma\) is connected, and passes through all vertices of \(W\), since the degree of each vertex is at most \(2\).
    Therefore, \(W = \{w_1, w_2, \dots, w_k\}\),
    where \(w_1, w_2, w_3, \dots, w_{k-1}, w_k\) are the vertices of a walk from \(w_1\) to \(w_k\) in \(\gamma\).

    For each \(i \in \{1, 2, \dots, k-1\}\), the undirected edge \(w_i w_{i+1}\) is an edge of \(\gamma\).
    Therefore, at least one of the directed edges \(w_i w_{i+1}\) and \(w_{i+1} w_i\) is an edge of \(h|_W\).
    If only one of \(w_i w_{i+1}\) and \(w_{i+1} w_i\) is an edge of \(h|_W\), say \(w_i w_{i+1}\), then there is a (simple) directed path from \(w_{i+1}\) to \(w_i\) in \(h|_W\).
    Since \(W\) is a strongly connected component of \(h\), there exists a (simple) cycle in \(h|_W\) of length at least \(3\).
    Since \(\graph\) is a directed subgraph of a bidirected tree, this is impossible.
    Therefore, both \(w_i w_{i+1}\) and \(w_{i+1} w_i\) are edges of \(h|_W\).

    Since this gives the total count of edges in \(h|_W\),
    these are the only edges of \(h|_W\), and \(h|_W\) is the multifunction along the path \(w_1, w_2, w_3, \dots, w_{k-1}, w_k\) from \(I\) to \(J\).
\end{proof}

Note that the edges of \(h|_W\) are uniquely determined by \(W\) and \(\gamma\).

\begin{prop}
    Let \(p\) be a multifunction along a path \(P = w_1w_2\dots w_k\) from \(I\) to \(J\).
    Then,
    \begin{enumerate}[label=\textup{(\roman*)}, leftmargin = .6in]
        \item \label{item:uniqueness-I-J}
        there exist unique acyclic functions \(p_I\) on \(W_I(P)\) and \(p_J\) on \(W_J(P)\) such that \(u = p_I + p_J\),
        \item \label{item:nonexistence-I-cup-J}
        there does not exist acyclic functions \(p_{I \cup J}\) on \(W_{I \cup J}(P)\) and \(p_{I \cap J}\) on \(W_{I \cap J}(P)\) such that \(u = p_{I \cup J} + p_{I \cap J}\).
    \end{enumerate}
\end{prop}

\begin{proof}[Proof (existence of \ref{item:uniqueness-I-J})]
    Define \(p_I\) on \(W_I(P)\) as the multifunction consisting of the edges \(w_1 w_2\), \(w_2 w_3\), \dots, \(w_{k-1} w_k\),
    and define \(p_J\) on \(W_J(P)\) as the multifunction consisting of the edges \(w_2 w_1\), \(w_3 w_2\), \dots, \(w_k w_{k-1}\).
    Then, \(p = p_I + p_J\).
\end{proof}

\begin{figure}[htbp]
    \centering
    \begin{tikzpicture}[x=1cm]
        \draw[fill=blue!40, draw=black, fill opacity=.5] (-.5, -.7) rectangle (5.5, .6);
        \draw[fill=red!40, draw=black, fill opacity=.5] (.5, -.6) rectangle (6.5, .7);
        \node[shape=circle, draw=black] (1) at (0,0){};
        \node[shape=circle, draw=black] (2) at (1,0){};
        \node[shape=circle, draw=black] (3) at (2,0){};
        \node[] (4) at (3,0){\small \(\cdots\)};
        \node[shape=circle, draw=black] (5) at (4,0){};
        \node[shape=circle, draw=black] (6) at (5,0){};
        \node[shape=circle, draw=black] (7) at (6,0){};
        \draw[->, blue!70!black] (1) edge[bend left] (2);
        \draw[->, blue!70!black] (2) edge[bend left] (3);
        \draw[->, blue!70!black] (3) edge[bend left] (4);
        \draw[->, blue!70!black] (4) edge[bend left] (5);
        \draw[->, blue!70!black] (5) edge[bend left] (6);
        \draw[->, blue!70!black] (6) edge[bend left] (7);
        \draw[->, red!70!black] (2) edge[bend left] (1);
        \draw[->, red!70!black] (3) edge[bend left] (2);
        \draw[->, red!70!black] (4) edge[bend left] (3);
        \draw[->, red!70!black] (6) edge[bend left] (5);
        \draw[->, red!70!black] (5) edge[bend left] (4);
        \draw[->, red!70!black] (7) edge[bend left] (6);
    \end{tikzpicture}
    \caption{Illustration of construction of \(p_I\) on \(W_I(P)\) (in blue) and \(p_J\) on \(W_J(P)\) (in red) for which \(p = f + g\).}
\end{figure}

\begin{proof}[Proof (uniqueness of \ref{item:uniqueness-I-J})]
    Let \(p_I\) and \(p_J\) be acyclic functions on \(W_I\) and \(W_J\) such that \(u = p_I + p_J\).
    Note that \(w_1 w_2 \in p_I\).
    If \(w_i w_{i+1} \in p_I\) for all \(i \in \{2, \dots, k-1\}\), then we obtain the construction above.
    Otherwise, let \(i \in \{2, \dots, k-1\}\) be the smallest index for which \(w_i w_{i+1} \notin p_I\).
    On the one hand, since \(w_i\) sends an edge in \(p_I\), we obtain \(w_i w_{i-1} \in p_I\).
    On the other hand, by minimality, we have \(w_{i - 1} w_i \in p_I\).
    Hence, \(p_I\) is not acyclic, a contradiction.
    Therefore, the only acyclic functions \(p_I\) and \(p_J\) on \(W_I\) and \(W_J\) such that \(p = p_I + p_J\) are the ones constructed above.
\end{proof}

\begin{proof}[Proof (nonexistence of \ref{item:nonexistence-I-cup-J})]
    Assume there exist acyclic functions \(p_{I \cup J}\) on \(W_{I \cup J}\) and \(p_{I \cap J}\) on \(W_{I \cap J}\) such that \(p = p_{I \cup J} + p_{I \cap J}\).
    Note that \(w_1 w_2 \in p_{I \cup J}\) and \(w_{k} w_{k-1} \in p_{I \cup J}\).
    Let \(i \in \{1, 2, \dots, k-1\}\) be the largest index such that \(w_i w_{i+1} \in p_{I \cup J}\).
    If \(i = k-1\), then \(w_{k-1} w_k, w_k w_{k-1} \in p_{I \cup J}\) so \(p_{I \cup J}\) is not acyclic, a contradiction.
    Assume \(i < k-1\).
    Then, \(w_i w_{i+1} \in p_{I \cup J}\).
    By maximality, \(w_{i+1} w_{i+2} \notin p_{I \cup J}\).
    Since \(w_{i+1}\) sends an edge in \(p_{I \cup J}\), we obtain \(w_{i+1} w_i \in p_{I \cup J}\).
    Hence, \(p_{I \cup J}\) is not acyclic, a contradiction.
    Therefore, there do not exist acyclic functions \(p_{I \cup J}\) and \(p_{I \cap J}\) on \(W_{I \cup J}\) and \(W_{I \cap J}\) such that \(p = p_{I \cup J} + p_{I \cap J}\).
\end{proof}

The previous proposition can be stacked to obtain analogous proposition for multifunctions along families of disjoint paths.

\begin{coro} \label{corollary:disjointpaths}
    Let \(p\) be the multifunction along a family \(\mathcal{P}\) of disjoint paths from \(I\) to \(J\).
    Then,
    \begin{enumerate}[label=\textup{(\roman*)}]
        \item there exist unique acyclic functions \(p_I\) on \(W_I(\mathcal{P})\) and \(p_J\) on \(W_J(\mathcal{P})\) such that \(p = p_I + p_J\), and
        \item there does not exist acyclic functions \(p_{I \cup J}\) on \(W_{I \cup J}(\mathcal{P})\) and \(p_{I \cap J}\) on \(W_{I \cap J}(\mathcal{P})\) such that \(p = p_{I \cup J} + p_{I \cap J}\).
    \end{enumerate}
\end{coro}

\subsubsection{Partition of \texorpdfstring{\(\mathcal{F}\)}{F} into \texorpdfstring{\(\mathcal{F}_{\mathcal{P}}\)}{FP}}

Given a multifunction \(p\) along a family of disjoint paths \(\mathcal{P}\) from \(I\) to \(J\),
define \(\mathcal{F}_{\mathcal{P}}\) to be the set of pairs of acyclic functions \((\bar{f}, \bar{g})\) such that
\begin{itemize}
    \item \(\bar{f}\) is an acyclic function on \(I \setminus W_I(\mathcal{P})\),
    \item \(\bar{g}\) is an acyclic function on \(J \setminus W_J(\mathcal{P})\), and
    \item \(\bar{f} + \bar{g}\) has no strongly connected components which are not entirely contained in \(I\) nor entirely contained in \(J\).
\end{itemize}

\begin{lemm} \label{lemma:partition-f-fp}
    \begin{equation*}
        \sum_{(f,g)\in\mathcal{F}} \normalizedweight(f+g)
        = \sum_{\mathcal{P}} \sum_{(\bar{f}, \bar{g}) \in \mathcal{F}_{\mathcal{P}}}\normalizedweight(\bar{f} + \bar{g}),
    \end{equation*}
    where the sum is over all families of disjoint paths \(\mathcal{P}\) from \(I\) to \(J\).
\end{lemm}

\begin{proof}
    It suffices to construct a \(\normalizedweight\)-preserving bijection between
    \begin{equation*}
        \mathcal{F} 
        \qquad \text{and} \qquad 
        \bigsqcup_{\mathcal{P}} \mathcal{F}_{\mathcal{P}}.
    \end{equation*}

    Define \(\alpha \colon \mathcal{F} \to \bigsqcup_{\mathcal{P}} \mathcal{F}_{\mathcal{P}}\) as follows.
    Given \((f, g) \in \mathcal{F}\),
    consider the strongly connected components of \(f+g\).
    By Proposition~\ref{proposition:restriction-strongly-connected},
    each strongly connected component that is not entirely contained in \(I\) or \(J\) has an associated path from \(I\) to \(J\).
    Let \(\mathcal{P}\) be the family of disjoint paths associated with these components,
    and let \(p\) be the multifunction along \(\mathcal{P}\).
    Recall the existence of unique acyclic functions \(p_I\) on \(W_I(\mathcal{P})\) and \(p_J\) on \(W_J(\mathcal{P})\) such that \(p = p_I + p_J\) from Corollary~\ref{corollary:disjointpaths}.
    Since the \(p \subset f+g\), we have \(p_I \subset f\) and \(p_J \subset g\).
    Define \(\alpha(f, g) = (f - p_I, g - p_J)\), which are acyclic functions on \(I \setminus W_I(\mathcal{P})\) and \(J \setminus W_J(\mathcal{P})\).
    The map is well-defined, since \(f - p_I + g - p_J = f + g - p\) has no strongly connected components which are not entirely contained in \(I\) nor entirely contained in \(J\).
    Moreover, since \(\normalizedweight(p) = 1\),
    \begin{equation*}
        \normalizedweight(f - p_I + g - p_J) = \normalizedweight(f + g - p) = \normalizedweight(f + g),
    \end{equation*}
    and consequently \(\alpha\) is \(\normalizedweight\)-preserving.

    Define \(\beta \colon \bigsqcup_{\mathcal{P}} \mathcal{F}_{\mathcal{P}} \to \mathcal{F}\) as follows.
    Let \((\bar{f}, \bar{g}) \in \mathcal{F}_{\mathcal{P}}\) for some family \(\mathcal{P}\) of disjoint paths from \(I\) to \(J\).
    Let \(p\) be the multifunction along \(\mathcal{P}\).
    Define \(\beta(\bar{f}, \bar{g}) = (\bar{f} + p_I, \bar{g} + p_J)\).
    Note that, since \(\normalizedweight(p) = 1\),
    \begin{equation*}
        \normalizedweight(\bar{f} + p_I + \bar{g} + p_J) = \normalizedweight(\bar{f} + \bar{g} + p) = \normalizedweight(\bar{f} + \bar{g}),
    \end{equation*}
    and consequently \(\beta\) is \(\normalizedweight\)-preserving.

    Direct computation shows that \(\alpha\) and \(\beta\) are inverses of each other.
    Therefore, \(\alpha\) is a \(\normalizedweight\)-preserving bijection between \(\mathcal{F}\) and \(\bigsqcup_{\mathcal{P}} \mathcal{F}_{\mathcal{P}}\), as desired.
\end{proof}

\subsubsection{Bijection between \texorpdfstring{\(\mathcal{F}_{\mathcal{P}}\) and \(\mathcal{H}_{\mathcal{P}}\)}{FP and HP}}

Given a multifunction along a family of disjoint paths \(\mathcal{P}\) from \(I\) to \(J\),
define \(\mathcal{H}_{\mathcal{P}}\) be the set of pairs \((f', g')\) such that
\(f'\) is an acyclic function on \((I \cup J) \setminus W_{I \cup J}(\mathcal{P})\), and
\(g'\) is an acyclic function on \((I \cap J) \setminus W_{I \cap J}(\mathcal{P})\).
From Corollary~\ref{corollary:disjointpaths},
\(f' + g'\) has no strongly connected components which are not entirely contained in \(I\) nor entirely contained in \(J\).
By definition, we have
\begin{equation*}
    Y_{I \cup J \setminus W_{I \cup J}(\mathcal{P})} Y_{I \cap J \setminus W_{I \cap J}(\mathcal{P})}
    =
    \sum_{(f', g') \in \mathcal{H}_{\mathcal{P}}} \normalizedweight(f' + g').
\end{equation*}

\begin{prop} \label{proposition:nwt-fp-hp}
    Let \(\mathcal{P}\) be a family of disjoint paths from \(I\) to \(J\).
    Then,
    \begin{equation*}
    \sum_{(f,g) \in \mathcal{F}_{\mathcal{P}}} \normalizedweight(f+g)
    =
    \sum_{(f',g') \in \mathcal{H}_{\mathcal{P}}} \normalizedweight(f'+g').
    \end{equation*}
\end{prop}

To prove this proposition, we define a \(\normalizedweight\)-preserving bijection between \(\mathcal{F}_{\mathcal{P}}\) and \(\mathcal{H}_{\mathcal{P}}\).

\begin{proof}[Proof of Proposition~\ref{proposition:nwt-fp-hp}]
We define the map \(\phi \colon \mathcal{F}_{\mathcal{P}} \to \mathcal{H}_{\mathcal{P}}\) as follows.
Given \((f, g) \in \mathcal{F}_{\mathcal{P}}\),
partition the vertex set \((I \cup J) \setminus W_{I \cup J}(\mathcal{P})\) into strongly connected components of \(f+g\).
Let \(V_I\) be the set of vertices in components entirely contained in \(I\), and let \(V_J\) be the set of vertices in components entirely contained in \(J\) but not in \(I\).
By definition of \(\mathcal{F}_{\mathcal{P}}\), there are no components of \(f+g\) which are not entirely contained in \(I\) nor entirely contained in \(J\).
Hence, \(V_I\) and \(V_J\) form a partition of \((I \cup J) \setminus W_{I \cup J}(\mathcal{P})\).
Define \(\phi(f, g) = (f', g')\) where
\begin{equation*}
    f' = f|_{V_I} + g|_{V_J} \qquad \text{and} \qquad g' = f|_{V_J} + g|_{V_I}.
\end{equation*}

\begin{claim} \label{claim:phi-well-defined}
    The map \(\phi\) is a well-defined \(\normalizedweight\)-preserving map from \(\mathcal{F}_{\mathcal{P}}\) to \(\mathcal{H}_{\mathcal{P}}\).
\end{claim}

\begin{proof}
Note that \(f' + g' = f + g\).
Hence, \(\phi\) is \(\normalizedweight\)-preserving.
We check that \(\phi\) is well-defined, that is, \((f', g') \in \mathcal{H}_{\mathcal{P}}\).
First, we check that \(f'\) and \(g'\) are functions on \((I \cup J) \setminus W_{I \cup J}(\mathcal{P})\) and on \((I \cap J) \setminus W_{I \cap J}(\mathcal{P})\), respectively.
Assume \(v \in (I \cap J) \setminus W_{I \cap J}(\mathcal{P})\).
If \(v \in V_I\), then \(v\) has exactly one outgoing edge in \(f'\) as part of \(f|_{V_I}\), 
and exactly one outgoing edge in \(g'\) as part of \(g|_{V_I}\).
If \(v \in V_J\), then \(v\) has exactly one outgoing edge in \(f'\) as part of \(g|_{V_J}\),
and exactly one outgoing edge in \(g'\) as part of \(f|_{V_J}\).
Assume \(v \in (I \cup J) \setminus (I \cap J) \setminus W_{I \cup J}(\mathcal{P})\).
If \(v \in V_I\), then \(v\) has exactly one outgoing edge in \(f'\) as part of \(f|_{V_I}\), and no outgoing edges in \(g'\).
If \(v \in V_J\), then \(v\) has exactly one outgoing edge in \(g'\) as part of \(g|_{V_J}\), and no outgoing edges in \(f'\).
Therefore, \(f'\) and \(g'\) are functions on \((I \cup J) \setminus W_{I \cup J}(\mathcal{P})\) and on \((I \cap J) \setminus W_{I \cap J}(\mathcal{P})\), respectively.

Now, we check that \(f'\) and \(g'\) are acyclic.
If \(c \subset f'\) is a cycle, then \(c \subset f' + g' = f + g\).
Since \(f+g\) has no strongly connected components which are not entirely contained in \(I\) nor entirely contained in \(J\), then \(c\) is entirely contained in \(I\) or \(J\).
If \(c\) is entirely contained in \(I\), then \(f' = f|_{V_I} + g|_{V_J}\) implies that \(c \subset f|_{V_I}\), a contradiction because \(f\) is acyclic.
By an analogous argument, if \(c\) is entirely contained in \(J\), then \(g' = f|_{V_J} + g|_{V_I}\) implies that \(c \subset g|_{V_I}\), a contradiction because \(g\) is acyclic.
Therefore, \(f'\) is an acyclic function on \((I \cup J) \setminus W_{I \cup J}(\mathcal{P})\).
By an analogous argument, \(g'\) is an acyclic function on \((I \cap J) \setminus W_{I \cap J}(\mathcal{P})\).
Therefore, \((f', g') \in \mathcal{H}_{\mathcal{P}}\).
\end{proof}

Moreover, define the map \(\psi \colon \mathcal{H}_{\mathcal{P}} \to \mathcal{F}_{\mathcal{P}}\) as follows.
Given \((f', g') \in \mathcal{H}_{\mathcal{P}}\),
partition the vertex set \((I \cup J) \setminus W_{I \cup J}(\mathcal{P})\) into strongly connected components of \(f'+g'\).
Let \(V_I\) be the set of vertices in components entirely contained in \(I \setminus W_I(\mathcal{P})\),
and let \(V_J\) be the set of vertices in components entirely contained in \(J \setminus W_J(\mathcal{P})\) but not in \(I \setminus W_I(\mathcal{P})\).
Recall that there are no components of \(f'+g'\) which are not entirely contained in \(I\) nor entirely contained in \(J\).
Hence, \(V_I\) and \(V_J\) form a partition of \((I \cup J) \setminus W_{I \cup J}(\mathcal{P})\).
Define \(\psi(f', g') = (f, g)\) where
\begin{equation*}
    f = f'|_{V_I} + g'|_{V_J} \qquad \text{and} \qquad g = f'|_{V_J} + g'|_{V_I}.
\end{equation*}

\begin{claim} \label{claim:psi-well-defined}
    The map \(\psi\) is a well-defined \(\normalizedweight\)-preserving map from \(\mathcal{H}_{\mathcal{P}}\) to \(\mathcal{F}_{\mathcal{P}}\).
\end{claim}

The proof of Claim~\ref{claim:psi-well-defined} is analogous to the proof of Claim~\ref{claim:phi-well-defined}.

\begin{claim}
    The maps \(\phi\) and \(\psi\) are inverses of each other.
\end{claim}
\begin{proof}
    Let \((f, g) \in \mathcal{F}_{\mathcal{P}}\), and let \((f', g') = \phi(f, g) \in \mathcal{H}_{\mathcal{P}}\).
    Since \(f+g = f'+g'\), the sets of strongly connected components of \(f+g\) and \(f'+g'\) are the same, i.e., \(V_I\) and \(V_J\) are the same in the definitions of \(\phi\) and \(\psi\).
    Therefore, if \((f'', g'') = \psi(f', g')\),
    then
    \begin{align*}
        f'' = f'|_{V_I} + g'|_{V_J} = (f|_{V_I} + g|_{V_J})|_{V_I} + (f|_{V_J} + g|_{V_I})|_{V_J} = f|_{V_I} + f|_{V_J} = f.
    \end{align*}
    Therefore, \(\psi \circ \phi = \operatorname{id}_{\mathcal{F}_{\mathcal{P}}}\).
    By an analogous argument, \(\phi \circ \psi = \operatorname{id}_{\mathcal{H}_{\mathcal{P}}}\).
    Therefore, \(\phi\) and \(\psi\) are inverses of each other.
\end{proof}

Finally, the map \(\phi\) is a \(\normalizedweight\)-preserving bijection between \(\mathcal{F}_{\mathcal{P}}\) and \(\mathcal{H}_{\mathcal{P}}\).
As a consequence, we have
\begin{equation*} \label{eq:sum-Fp-equals-sum-Hp}
    \sum_{(f,g) \in \mathcal{F}_{\mathcal{P}}} \normalizedweight(f+g) = \sum_{(f',g') \in \mathcal{H}_{\mathcal{P}}} \normalizedweight(f'+g'). \qedhere
\end{equation*}
\end{proof}

\subsubsection{Aggregating over all families of disjoint paths}

We are finally ready to prove Proposition~\ref{proposition:expansionfortrees}.

\begin{prop*}[Proposition~\ref{proposition:expansionfortrees}, repeated]
    Let \(I, J \subset \vertices\).
    The expansion of \(Y_IY_J\) into cluster monomials is given by
    \begin{equation*}
        Y_IY_J = \sum_{\mathcal{P}} Y_{I \cup J \setminus W_{I \cup J}(\mathcal{P})}Y_{I \cap J \setminus W_{I \cap J}(\mathcal{P})},
    \end{equation*} 
    where the sum is taken over all families \(\mathcal{P}\) of disjoint paths from \(I\) to \(J\).
\end{prop*}

\begin{proof}[Proof of Proposition~\ref{proposition:expansionfortrees}]
    Proposition~\ref{proposition:nwt-fp-hp} implies that, for each family \(\mathcal{P}\) of disjoint paths from \(I\) to \(J\),
    \begin{equation*}
        \sum_{(f,g) \in \mathcal{F}_{\mathcal{P}}} \normalizedweight(f+g)
        =
        Y_{I \cup J \setminus W_{I \cup J}(\mathcal{P})}
        Y_{I \cap J \setminus W_{I \cap J}(\mathcal{P})}.
    \end{equation*}
    Adding over all \(P \in \mathcal{P}\) and using Lemma~\ref{lemma:partition-f-fp}, we obtain
    \begin{equation*}
        Y_IY_J
        =
        \sum_{(f,g) \in \mathcal{F}} \normalizedweight(f+g)
        =
        \sum_{\mathcal{P}} \sum_{(f,g) \in \mathcal{F}_{\mathcal{P}}} \normalizedweight(f+g)
        =
        \sum_{\mathcal{P}} Y_{I \cup J \setminus W_{I \cup J}(\mathcal{P})}Y_{I \cap J \setminus W_{I \cap J}(\mathcal{P})}. \qedhere
    \end{equation*}
\end{proof}

\subsection{Nonnegative expansion of \texorpdfstring{\(Y\)}{Y}-monomials}\label{subsection:nonnegativeexpansionofY}

In this subsection, we show that cluster \(Y\)-monomials span \(Y\)-monomials with nonnegative coefficients when \(\graph\) is a bidirected tree, as stated in Proposition~\ref{proposition:tree}.
This is the analogous result to Theorem~\ref{theorem:Yclusters-span} in the context of nonnegativity for bidirected trees.

\begin{defi}[Dominance order]
    Let \(\mathcal{A}\) and \(\mathcal{B}\) be two multisets whose elements are positive integers,
    and have equal sum.
    We say that \(\mathcal{A}\) \vocab{dominates} \(\mathcal{B}\),
    denoted \(\mathcal{A} \trianglerighteq \mathcal{B}\), if for all \(k \in \positiveintegers\),
    the sum of the \(k\) largest elements of \(\mathcal{A}\) is
    at least the sum of the \(k\) largest elements of \(\mathcal{B}\).
    (If \(k\) is larger than the size of a multiset, then the sum of the \(k\) largest elements is the sum of all elements.)
\end{defi}

For example, \(\{4, 4, 1\} \trianglerighteq \{4, 3, 1, 1\}\), since \(4 \geq 4\), \(4 + 4 \geq 4 + 3\), \(4 + 4 + 1 = 9 \geq 4 + 3 + 1\), and \(4 + 4 + 1 \geq 4 + 3 + 1 + 1 \).
It is known that the dominance order is a partial order on multisets whose elements are positive integers and have a fixed sum \(n\).

Dominance order is useful because,
given a multicollection \(\mathcal{S}\)
with sets \(I\) and \(J\) of vertices,
if \(\mathcal{T}\) is the nested multicollection obtained by replacing \(I\) and \(J\) by \(I \cup J\) and \(I \cap J\),
then the multiset of sizes of sets in \(\mathcal{T}\) is dominated by the multiset of sizes of sets in \(\mathcal{S}\), that is,
\( \{|T| : T \in \mathcal{T}\} \trianglerighteq \{|S| : S \in \mathcal{S}\} \),
where equality holds if and only if \(\{I, J\} = \{I \cup J, I \cap J\}\),
that is, \(I \subset J\) or \(J \subset I\).

Finally, note that both the dominance order and the reverse dominance order are well-founded on multisets of positive integers with a fixed sum, meaning that any non-empty set of multisets with a fixed sum has a minimal element with respect to the dominance order and a maximal element with respect to the reverse dominance order.
This means that we can apply induction on the dominance order and the reverse dominance order on multisets of positive integers with a fixed sum.

\begin{prop} \label{proposition:tree}
    Assume \(\graph\) is a bidirected tree.
    Then, each \(Y\)-monomial of \(\graphLPalgebra{\graph}\) is a nonnegative integer linear combination of cluster \(Y\)-monomials.

    More precisely, the \(Y\)-monomial \(\monomial{\varnothing}{\mathcal{S}}\) is a nonnegative integer linear combination of cluster \(Y\)-monomials \(\monomial{\varnothing}{\mathcal{N}}\), where \(\mathcal{N}\) is a nested multicollection such that the vertices in sets of \(\mathcal{N}\) are a subset of the vertices in sets of \(\mathcal{S}\).
\end{prop}

\begin{proof}[Proof of Proposition~\ref{proposition:tree}]
    We apply strong induction on \(\left|\sum_{S \in \mathcal{S}} S\right|\), the number of vertices in a set of \(\mathcal{S}\), counting multiplicities;
    and then we apply induction on \(\left\{|S| : S \in \mathcal{S}\right\}\), the multiset of sizes of sets in \(\mathcal{S}\), with respect to reverse dominance order.

    Let \(\mathcal{S}\) be a multicollection.
            Assume, by induction hypothesis, that the result holds for all nested multicollections \(\mathcal{R}\) such that \(\left|\sum_{R \in \mathcal{R}} R\right| < \left|\sum_{S \in \mathcal{S}} S\right|\).
    Also assume, by induction hypothesis, that the result holds for all nested multicollections \(\mathcal{R}\) such that \(\left|\sum_{R \in \mathcal{R}} R\right| = \left|\sum_{S \in \mathcal{S}} S\right|\) and \(\mathcal{R} \triangleright \mathcal{S}\).

                        Recall Lemma~\ref{lemma:nestedcollection:bidirected},
    which states that a multicollection \(\mathcal{N}\) is nested if and only if it satisfies conditions \ref{item:nestedcollection:disjoint}, \ref{item:nestedcollection:stronglyconnected:bidirected}, and \ref{item:nestedcollection:connected}.
    We split into whether \(\mathcal{S}\) satisfies conditions \ref{item:nestedcollection:disjoint} and \ref{item:nestedcollection:stronglyconnected:bidirected}.

    If \ref{item:nestedcollection:disjoint} or \ref{item:nestedcollection:stronglyconnected:bidirected} fails,
    then there exist \(I, J \in \mathcal{S}\)
    such that \(I \not\subset J\) and \(J \not\subset I\).
    This implies that \(\{I \cup J, I \cap J\} \neq \{I, J\}\).
    Given a family \(\mathcal{P}\) of disjoint paths from \(I\) to \(J\),
    let \(\mathcal{R}_{\mathcal{P}}\) denote the nested multicollection obtained by replacing \(I\) and \(J\) by \((I \cup J) \setminus W_{I \cup J}(\mathcal{P})\) and \((I \cap J) \setminus W_{I \cap J}(\mathcal{P})\).
    By Proposition~\ref{proposition:expansionfortrees},
    \begin{equation*}
        \monomial{\varnothing}{\mathcal{S}}
        =
        \sum_{\mathcal{P}} \monomial{\varnothing}{\mathcal{R}_{\mathcal{P}}}.
    \end{equation*}
    For each non-empty \(\mathcal{P}\), the nested multicollection \(\mathcal{R}_{\mathcal{P}}\) has fewer vertices than \(\mathcal{S}\),
    that is, \(\left|\sum_{R \in \mathcal{R}_{\mathcal{P}}} R\right| < \left|\sum_{S \in \mathcal{S}} S\right|\).
    Therefore, by the induction hypothesis, each \(\monomial{\varnothing}{\mathcal{R}_{\mathcal{P}}}\) is a nonnegative integer linear combination of cluster \(Y\)-monomials.
    For \(\mathcal{P} = \varnothing\), the nested multicollection \(\mathcal{R}_{\varnothing}\) has the same number of vertices as \(\mathcal{S}\),
    that is, \(\left|\sum_{R \in \mathcal{R}_{\varnothing}} R\right| = \left|\sum_{S \in \mathcal{S}} S\right|\),
    but additionally it satisfies \(\mathcal{R}_{\varnothing} \triangleright \mathcal{S}\).
    Therefore, by the induction hypothesis, \(\monomial{\varnothing}{\mathcal{R}_{\varnothing}}\) is a nonnegative integer linear combination of cluster \(Y\)-monomials.
    Hence, if follows that \(\monomial{\varnothing}{\mathcal{S}}\) is also a nonnegative integer linear combination of cluster \(Y\)-monomials.

    Otherwise, assume that conditions \ref{item:nestedcollection:disjoint} and \ref{item:nestedcollection:stronglyconnected:bidirected} hold.
    Let \(\mathcal{N}\) be the multicollection
    obtained from \(\mathcal{S}\) by replacing each set \(I \in \mathcal{S}\) with the connected components of \(I\).
    Now, \(\mathcal{N}\) satisfies \ref{item:nestedcollection:connected},
    while \ref{item:nestedcollection:disjoint} and \ref{item:nestedcollection:stronglyconnected:bidirected} are preserved from \(\mathcal{S}\).
    Therefore,
    by Lemma~\ref{lemma:nestedcollection:bidirected},
    \(\mathcal{N}\) is a nested multicollection.
    Moreover, by Lemma~\ref{lemma:y:product}, we have \(\monomial{\varnothing}{\mathcal{S}} = \monomial{\varnothing}{\mathcal{N}}\),
    so \(\monomial{\varnothing}{\mathcal{S}}\) equals to a cluster \(Y\)-monomial, as desired.
\end{proof}

\subsection{Nonnegative expansion of monomials}

In this subsection,
we prove Theorem~\ref{theorem:tree} as a corollary of Proposition~\ref{proposition:tree}.
This is very similar to the proof of Theorem~\ref{theorem:spanningset} as a corollary of Theorem~\ref{theorem:Yclusters-span},
with the additional care of making sure that the coefficients are nonnegative.

\begin{theo*}[Theorem~\ref{theorem:tree}, repeated]
    Assume that \(\graph\) is a bidirected tree.
    Then, each monomial of \(\graphLPalgebra{\graph}\) is a linear combination of cluster monomials with over \(\coefficientring\) with nonnegative coefficients.
\end{theo*}

The proof of Theorem~\ref{theorem:tree} is completely analogous to the proof of Theorem~\ref{theorem:spanningset}:
replace all instances of ``\emph{integer linear combination}'' to ``\emph{nonnegative integer linear combination}'',
replace all instances of ``\emph{linear combination over \(\coefficientring\)}'' to ``\emph{linear combination over \(\coefficientring\) with nonnegative coefficients}'',
and replace the application of Theorem~\ref{theorem:Yclusters-span}
to an application of Proposition~\ref{proposition:tree}.

\subsubsection*{Acknowledgements}

We thank the project's mentor Pasha Pylyavskyy and the project's teaching assistant Robbie Angarone for their guidance and support.
We also thank an anonymous referee for their helpful comments.

\printbibliography

\setlength{\parindent}{0em}

\end{document}